\documentstyle[12pt]{article}
\newfam\msbfam
\font\tenmsb=msbm10    \textfont\msbfam=\tenmsb \font\sevenmsb=msbm7
\scriptfont\msbfam=\sevenmsb \font\fivemsb=msbm5
\scriptscriptfont\msbfam=\fivemsb
\def\Bbb{\fam\msbfam \tenmsb}

\def\rr{{\Bbb R}}
\def\rz{{{\rr}^n}}
\def\zz{{\Bbb Z}}
\def\nn{{\Bbb N}}
\def\cc{{\Bbb C}}

\def\fz{\infty}
\def\az{\alpha}
\def\supp{{\rm{\ supp\ }}}
\def\loc{{\rm{\ loc\ }}}

\def\ez{\epsilon}
\def\bz{\beta}
\def\pz{\partial}
\def\tz{\theta}
\def\sz{\sigma}
\def\vz{\varphi}
\def\dz{\delta}
\def\gz{\gamma}

\def\lz{\lambda}
\def\oz{\Omega}
\def\supp{{\rm supp}}
\def\loc{{\rm loc}}
\def\wt{\widetilde}
\def\wz{\omega}
\def\l{\left}
\def\r{\right}

\def\dsum{\displaystyle\sum}

\def\dint{\displaystyle\int}
\def\dfrac{\displaystyle\frac}
\def\dsup{\displaystyle\sup}

\def\dinf{\displaystyle\inf}

\newtheorem{thm}{\hskip\parindent Theorem}

\newtheorem{lem}{\hskip\parindent Lemma}

\newtheorem{prop}{\hskip\parindent Proposition}

\newtheorem{cor}{\hskip\parindent Corollary}

\begin{document}

\baselineskip=15pt
\renewcommand{\arraystretch}{2}
\arraycolsep=1.2pt

\title{ Weighted
local Hardy spaces and their applications} {\footnotetext{
\hspace{-0.65
cm} 2000 Mathematics Subject  Classification: 42B20, 42B25.\\
The  research was supported  by the NNSF (10971002)and (10871003) of China.\\}

\author{ Lin Tang}
\date{}
\maketitle

{\bf Abstract}\quad  In this paper, we  study weighted local Hardy
spaces $h^p_\wz(\rz)$ associated with local weights which include
the classical Muckenhoupt weights. This setting includes the
classical local Hardy space theory of Goldberg \cite{g}, and the weighted
Hardy spaces of Bui \cite{bu}.

\bigskip

\begin{center}{\bf 1. Introduction }\end{center}
The theory of local Hardy space plays an important role in various
fields of analysis and partial differential equations; see
\cite{s, st1, ta,v}. In particular, pseudo-differential operators
are bounded on local Hardy spaces $h^p$ for $0<p\le 1$, but they
are not bounded on Hardy spaces $H^p$ for $0<p\le 1$; see
\cite{g}.

On the other hand, Bui \cite{b} studied the weighted version
$h_w^p$ of the local Hardy space $h^p$ considered by Goldberg
\cite{g}, where the weight $\wz$ is assumed to satisfy the
condition $(A_\fz)$ of Muckenhoupt. Recently, Rychkov
\cite{v} introduced and studied some properties of the weighted
Besov-Lipschitz and Triebel-Lizorkin spaces  with weights that are
locally in $A_p$ but may grow or decrease exponentially. Recently, Rychkov
\cite{v} studied the class of Triebel-Lizorkin $F^s_{p,q}$ spaces, which includes Hardy spaces as its
part. In fact, Rychkov explicitly identifies weighted local Hardy space $h^p_\wz$ with $F^0_{p,2}(\wz)$ in
Theorem 2.25 of \cite{v}. In particular, Rychkov \cite{v} extended a part of theory of
$A_\fz$-weighted local Hardy spaces developed in Bui \cite{bu} to
the $A_\fz^{loc}$ weights, where $A_\fz^{loc}$ weights denote
local $A_\fz$-weights which are non-doubling weights, and the
$A_\fz^{loc}$ weights include the $A_\fz$-weights.

The main purpose of this paper is twofold. The first goal is to
establish weighted atomic decomposition characterizations of
weighted local Hardy space $h_\wz^p$ with local weights. The second
goal is to show that  strong singular integrals and
Pseudodifferential operators and their commutators are bounded on
weighted local Hardy spaces.

The paper is organized as follows. In Section 2, we first recall
some notation and definitions concerning local weights and grand
maximal function; and we then obtain a basic approximation of the
identity result and the grand maximal function characterization
for $L^q_\wz$ with $q\in (q_\wz,\fz]$, where $q_\wz$ is the
critical of $\wz$. In Section 3, we introduce weighted local Hardy
spaces $h^p_{\wz,N}$ via grand maximal functions and weighted
atomic local Hardy spaces $h_{\wz}^{p,q,s}(\rz)$ for any
admissible triplet $(p,q,s)_\wz$, and study some properties of
these spaces. In Section 4, we establish the Calder\'on-Zygmund
decomposition associated with the grand maximal function. In
Section 5, we prove that for any admissible  triplet
$(p,q,s)_\wz$, $h^p_{\wz, N}(\rz)=h_{\wz}^{p,q,s}(\rz)$ with
equivalent norms.  Moreover, we prove that
$\|\cdot\|_{h_{\wz,fin}^{p,q,s}(\rz)}$ and
$\|\cdot\|_{H_\wz^p(\rz)}$ are equivalent quasi-norms on
$h^{p,q,s}_{\wz,fin}(\rz)$ with $q<\fz$, and we obtain criterions
for boundedness of sublinear operators in $h^p_\wz$ in Section 6.
Finally, in Section 7, we show that strong singular integrals and
Pseudodifferential operators and their commutators are bounded on
weighted local Hardy spaces by using weighted atomic
decompositions.

It is worth pointing out that  we can not adapt the methods in
\cite {bu} and \cite{g}, if $\wz$ is a local weight.
In fact,  adapting the same idea of (global)weighted Hardy spaces(
\cite{b,blyz,st,str}), we give a direct proof for
weighted atomic decompositions of weighted
 local Hardy spaces.

Throughout this paper, $C$ denotes the constants that are
independent of the main parameters involved but whose value may
differ from line to line. Denote by $\nn$ the set $\{1,2,\cdots\}$
and by $\nn_0$ the set $\nn\cup\{0\}$. By $A\sim B$, we mean that
there exists a constant $C>1$ such that $1/C\le A/B\le C$.

\begin{center}{\bf 2. Preliminaries  }\end{center}
We first introduce  weight classes $A_p^{loc}$ from \cite{v}.

Let $Q$ run through all cubes in $\rz$ (here and below only cubes
with sides parallel to the coordinate axes are considered), and let
$|Q|$ denote the volume of $Q$. We define the weight class
$A_p^{loc}(1<p<\fz$) to consists of all nonnegative locally integral
functions $\wz$ on $\rz$ for which
$$A_p^{loc}(\wz)=\dsup_{|Q|\le 1}\dfrac
1{|Q|^p}\dint_Q\wz(x)dx\l(\dint_Q\wz^{-p'/p}(x)dx\r)^{p/p'}<\fz ,
\  1/p+1/p'=1.\eqno(2.1)$$ The function $\wz$ is said to belong to
the weight class of $A_1^{loc}$ on $\rz$ for which
$$A_1^{loc}(\wz)=\dsup_{|Q|\le 1}\dfrac
1{|Q|}\dint_Q\wz(x)dx\l(\dsup_{y\in
Q}[\wz(y)]^{-1}\r)<\fz.\eqno(2.2)$$ {\bf Remark}: For any $C>0$ we
could have replaced $|Q|\le 1$ by $|Q|\le C$ in (2.1) and (2.2).

 In what follows,
$Q(x,t)$ denotes the cube centered at $x$ and of the sidelength $t$.
Similarly, given $Q=Q(x,t)$ and $\lz>0$, we will write $\lz Q$ for
the $\lz$-dilate cube, which is the cube with the same center $x$
and with sidelength $\lz t$. Given a Lebesgue measurable set $E$ and
a weight $\wz$,  let $\wz(E)=\int_E\wz dx$. For any $\wz\in
A_\fz^{loc}$, $L^p_\wz$ with $p\in(0,\fz)$ denotes the set of all
measurable functions $f$ such that
$$\|f\|_{L^p_\wz}\equiv \l(\dint_\rz|f(x)|^p\wz(x)dx\r)^{1/p}<\fz$$
and $L^\fz_\wz=L^\fz$.  We define the local Hardy-Littlewood
maximal operator by
$$M^{loc}f(x)=\dsup_{x\in Q:|Q|<1}\dfrac 1{|Q|}\dint_Q|f(y)|dy.$$
Similar to the
classical $A_p$ Muckenhoupt weights, we give some properties for
weights $\wz\in A^\loc_\fz:=\bigcup_{1\le p<\fz} A^{loc}_p$.
\begin{lem}\label{l2.1.}\hspace{-0.1cm}{\rm\bf 2.1.}\quad
Let $1\le p<\fz$, $\wz\in A_p^{loc}$, and $Q$ be a  unit cube,
i.e. $|Q|=1$. Then there exists a  $\bar\wz\in A_p$ so that
$\bar\wz=\wz$ on $Q$ and
\begin{enumerate}
\item[(i)]$A_p(\bar \wz)\le CA_p^{loc}(\wz).$
\item[(ii)]if $\wz\in A_p^{\loc}$, then there exists $\ez>0$ such that $\wz\in
A_{p-\ez}^{loc}(\wz)$ for $p>1$.
\item[(iii)]If $ 1\le p_1<p_2<\fz$, then $A_{p_1}^{loc}\subset
A_{p_2}^{loc}$. \item[(iv)] $\wz\in A_p^{loc}$ if and only if
$\wz^{-\frac 1{p-1}}\in A_{p'}^{loc}$.
\item[(v)] If $\wz\in A_p^{loc}$ for $1\le p<\fz$, then
$$\wz(tQ)\le exp(c_\wz t)\wz(Q)\quad (t\ge 1, |Q|=1).$$
\item[(vi)] the local Hardy-Littlewood maximal operator $M^{loc}$
is bounded on $L^p_\wz$ if $\wz\in A_p^{loc}$ with $p\in (1,\fz)$.
\item[(vii)]$M^{loc}$ is bounded from $L^1_\wz$ to $L^{1,\fz}_\wz$
if $\wz\in A_1^{loc}$.
\end{enumerate}
\end{lem}
Proof:  (i)-(vi) have been proved in \cite{v}. (vii) can be proved
by the standard method.

We remark that Lemma 2.1 is also true for $|Q|>1$ with $c$ depending
now on the size of $Q$. In addition, it is easy to see that
$A_p\subset A_p^{loc}$ for $p\ge 1$ and $e^{c|x|},\
(1+|x|\ln^\az(2+|x|))^\bz\in A_1^{\loc}$ with $\az\ge 0,\bz\in\rr$
and $c\in\rr$.

As a consequent of Lemma 2.1, we have following result.
\begin{cor}\label{c2.1.}\hspace{-0.1cm}{\rm\bf 2.1.}\quad
If $\wz\in A_\fz^{loc}$, then there exists a constant $C>0$ such that
$$\wz(2Q)\le \wz(Q)$$ if $|Q|<1$, and
$$\wz(Q(x_0,r+1))\le C\wz(Q(x_0,r))$$ if $|Q(x_0,r)|\ge 1$.
\end{cor}

From Lemma 2.1, for any given $\wz\in A^{\loc}_p$, define the
critical index of $\wz$ by
$$q_\wz\equiv\dinf\{p\in[1,\fz): \wz\in A_p^{loc}\}.\eqno(2.3)$$
Obviously, $q_\wz\in [1,\fz)$. If $q_\wz\in(1,\fz)$, then
$\wz\not\in A_{q_\wz}^{loc}$.

The symbols ${\cal D}(\rz)=C_0^\fz(\rz), {\cal D}'(\rz)$ is the
dual space of ${\cal D}(\rz)$. The multi-index notation is usual:
for $\az=(\az_1,\cdots,\az_n)$ and
$\pz^\az=(\pz/\pz_{x_1})^{\az_1}\cdots(\pz/\pz_{x_n})^{\az_n}$.
\begin{lem}\label{l2.2.}\hspace{-0.1cm}{\rm\bf 2.2.}\quad
Let $\wz\in A_\fz^{loc}, q_\wz$ be as in (2.3), and $p\in
(q_\wz,\fz]$. Then
\begin{enumerate}
\item[(i)] if $1/p+1/p'=1$, then ${\cal D}(\rz)\subset
L^{p'}_{\wz^{1/p-1}}(\rz)$; \item[(ii)] $L^p_\wz(\rz)\subset {\cal
D}'(\rz)$ and the inclusion is continuous.
\end{enumerate}
\end{lem}
Proof. We only prove the case $p<\fz$. The proof for the case
$p=\fz$ is easier and we omit the details. Since
$p\in(q_\wz,\fz)$, then $\wz\in A_p^{loc}$. Therefore, by the
definition of $A_p^{loc}$, for all ball $B=B(0,r)$ with radius $r$
and centered at $0$, we have
$$\dint_B[\wz(x)]^{-1/(p-1)}dx\le C[\wz(B)]^{-1/(p-1)}|B|^{p'}<\fz.$$
From this, for any $\vz\in {\cal D}(\rz)$ and $\supp\ \vz\subset
B$, we obtain
$$
\|\vz\|_{L^{p'}_{\wz^{-1/(p-1)}}(\rz)}\le C
\dint_B[\wz(x)]^{-1/(p-1)}dx<\fz.\eqno(2.4)
$$
For (ii), if $f\in L^p_\wz(\rz)$ and $\vz\in {\cal D}(\rz)$, by
H\"older inequality and (2.4), we have
$$|<f,\vz>|\le
\|f\|_{L^p_\wz(\rz)}\l(\dint_\rz|\vz(x)|^{p'}[\wz(x)]^{-1/(p-1)}dx\r)^{1/p'}\le
C\|f\|_{L^p_\wz(\rz)}.
$$ Thus, Lemma 2.2 is proved.

For $\vz\in {\cal D}(\rz)$ and $t>0$, set
$$\vz_t(x)=t^{-n}\vz\l(\frac xt\r).$$ It is easy to see that we have
the following results.
\begin{prop}\label{p2.1.}\hspace{-0.1cm}{\rm\bf 2.1.}\quad
Let $\vz\in {\cal D}(\rz)$ and $\dint_\rz \vz(x)dx=1$.
\begin{enumerate}
\item[(i)] For any $\Phi\in {\cal D}(\rz)$ and $f\in {\cal
D}'(\rz)$, $\Phi*\vz_t\to\Phi$ in ${\cal D}(\rz)$ as $t\to0$ and
$f*\vz_t\to f$ in ${\cal D}'(\rz)$ as $t\to0$. \item[(ii)] Let
$\wz\in A_\fz^{loc}$ and $q_\wz$ be as in (2.3). If $q\in
(q_\wz,\fz)$, then for any $f\in L^q_\wz(\rz)$, $f*\vz_t\to f$ in
$L^q_\wz(\rz)$ as $t\to0$.
\end{enumerate}
\end{prop}
Let $N\in \nn_0$ and
$$\begin{array}{cl}
{\cal M}_N^0 f(x)&=\dsup\{|\vz_t*f(x)|:\ 0<t<1,\vz\in{\cal D}(\rz),\ \int\vz\not=0,\\
&\qquad\supp\vz\subset B(0,1), \|D^\az\vz\|_\fz\le 1\ {\rm}\
|\az|\le N+1\}.\end{array}$$
$$\begin{array}{cl}
\bar{\cal M}_N^0 f(x)&=\dsup\{|\vz_t*f(x)|:\ 0<t<1,\vz\in{\cal D}(\rz),\ \int\vz\not=0,\\
&\qquad\supp\vz\subset B(0,2^{3(10+n)}), \|D^\az\vz\|_\fz\le 1\
{\rm}\ |\az|\le N+!\}.\end{array}$$ and

$$\begin{array}{cl}
{\cal  M}_N f(x)&=\dsup\{|\vz_t*f(z)|:\ |z-x|<t<1,\vz\in{\cal D}(\rz),\ \int\vz\not=0,\\
&\qquad\supp\vz\subset B(0,2^{3(10+n)})), \|D^\az\vz\|_\fz\le 1\
{\rm}\ |\az|\le N+1\}.\end{array}$$ For any $N\in \nn_0$ and
$x\in\rz$, obviously,
$${\cal M}_N^0 f(x)\le \bar{\cal M}_N^0 f(x)\le {\cal M}_N f(x).$$
For convenience, we write
$${\cal D}^0_N=\{\vz\in {\cal D}:\ \supp\vz\subset B(0,1),\ \int\vz\not=0,\
 \|D^\az\vz\|_\fz\le 1\ {\rm}\ |\az|\le N+1\},$$
and
$${\cal D}_N=\{\vz\in {\cal D}:\ \supp\vz\subset B(0,2^{3(10+n)}),\ \int\vz\not=0,\
 \|D^\az\vz\|_\fz\le 1\ {\rm}\ |\az|\le N+1\}.$$
\begin{prop}\label{p2.2.}\hspace{-0.1cm}{\rm\bf 2.2.}\quad
Let $N\ge 2$. Then
\begin{enumerate}
\item[(i)] There exists a positive $C$ such that for all $f\in
(L_{loc}^1(\rz)\bigcap{\cal D}'(\rz))$ and almost everywhere
$x\in\rz$, $|f(x)|\le {\cal  M}_N^0f(x)\le CM^{loc} f(x).$
\item[(ii)] If $\wz\in A_p^{loc}$ with $p\in (1,\fz)$, then $f\in
L^p_\wz(\rz)$ if and only if $f\in {\cal D}'(\rz)$ and ${\cal
M}_N^0f\in L^p_\wz$; moreover, $\|f\|_{L^p_\wz}\sim \|{\cal
M}_N^0f\|_{L^p_\wz}$. \item[(iii)] If $\wz\in A_1^{loc}$, then
${\cal M}_N^0$ is bounded from $L^1_\wz(\rz)$ to
$L^{1,\fz}_\wz(\rz)$.
\end{enumerate}
\end{prop}
The proof of (i) and (iii) is obvious, (ii) has been proved in
\cite{v},  we omit the details here.

\begin{center}{\bf 3. The grand maximal function definition of Hardy spaces  }\end{center}
In this section, we introduce  weighted local Hardy spaces via grand
maximal functions and weighted local Hardy spaces. Moreover, we
study some properties of these spaces.

Let $p\in (0,1]$, $\wz\in A_\fz^{loc}$, and $q_\wz$ be as in (2.3).
Set $$N_{p,\wz}=\max\{0,[n(\frac {q_\wz}p-1)]\}+2.$$ For each $N\ge
N_{p,\wz}$, the weighted local Hardy space is defined by
$$h_{\wz,N}^p(\rz)\equiv\l\{f\in {\cal D}'(\rz): {\cal M}^0_N(f)\in
L_\wz^p(\rz)\r\}.$$ Moreover, we define $\|f\|_{h^p_{\wz,
N}(\rz)}\equiv\|{\cal M}_N^0(f)\|_{L^p_\wz(\rz)}$. From  Theorem
2.24 in \cite{v}, we know that $\|{\cal
M}^0_N(f)\|_{L^p_\wz(\rz)}\sim\|\bar{\cal
M}_N^0(f)\|_{L^p_\wz(\rz)}\sim\|{\cal M}_N(f)\|_{L^p_\wz(\rz)}$.

For any integer $N,\bar N$ with $N_{p,\wz}\le N\le \bar N$, we have
$$h_{\wz,N_{p,\wz}}^p(\rz)\subset h_{\wz,N}^p(\rz)
\subset h_{\wz,\bar N}^p(\rz)$$ and the inclusions are continuous.

Notice that if $p\in (q_\wz,\fz]$ and $N\ge N_{p,\wz}=2$, then by
Proposition 2.2 (ii), we have  $h^p_{\wz,N}(\rz)=L^p_\wz(\rz)$
with equivalent norms. However, if $p\in (1,q_\wz)$, the element
of $h^p_{\wz,N}(\rz)$ may be a distribution, and hence,
$h^p_{\wz,N}(\rz)\not= L_\wz^p(\rz)$. But,
$(h^p_{\wz,N}(\rz))\bigcap L_{loc}^1(\rz))\subset L^p_\wz(\rz)$.
For applications considered in this paper, we concentrate only on
$h^p_{\wz,N}(\rz)$ with $p\in (0,1]$.

We introduce the following weighted atoms.

Let $\wz\in A_\fz^{loc}$ and $q_\wz$ be as in (2.3). A triplet
$(p,q,s)_\wz$ is called to be admissible, if $p\in (0,1]$,
$q\in(q_\wz,\fz]$ and $s\in\nn$ with $s\ge [n(\frac{q_\wz}p-1)]$.
A function $a$ on $\rz$ is said to be a $(p,q,s)_\wz-atom$ if
\begin{enumerate}
\item[(i)] $\supp \ a\subset Q$, \item[(ii)]
$\|a\|_{L^q_\wz(\rz)}\le [\wz( Q)]^{1/q-1/p}$.
\item[(iii)]$\dint_\rz a(x)x^\az dx=0$ for $\az\in (\nn_0)^n$ with
$|\az|\le s$, if $|Q|<1$.
\end{enumerate}
Moreover, we call $a$ is a $(p,q)_\wz$ single atom if
$\|a\|_{L^q_\wz(\rz)}\le [\wz(\rz)]^{1/q-1/p}$.

Let $\wz\in A_\fz^{loc}$ and $(p,q,s)_\wz$ be an admissible
triplet.  The weighted atomic local Hardy space
$h^{p,q,s}_{\wz}(\rz)$ is defined to be the set of all $f\in {\cal
D}'(\rz)$ satisfying that $f=\sum_{i=0}^\fz\lz_ia_i$ in ${\cal
D}'(\rz)$, where $\{\lz_i\}_{i\in\nn_0}\subset \cc, \
\sum_{i=0}^\fz |\lz_i|^p<\fz$ and $\{a_i\}_{i\in\nn}$ are
$(p,q,s)_\wz$-atom and $a_0$ is a $(p,q)_\wz$ single atom.
Moreover, the quasi-norm of $f\in h^{p,q,s}_{\wz}(\rz)$ is defined
by
$$\|f\|_{h^{p,q,s}_{\wz}(\rz)}\equiv\dinf\l\{\l[\dsum_{i=0}^\fz|\lz_i|^p\r]^{1/p}\r\},$$
where the infimum is taken over all the decompositions of $f$ as
above.

It is easy to see that if the triplets $(p,q,s)_\wz$ and $(p,\bar
q,\bar s)_\wz$ are admissible and  satisfy $\bar q\le q$ and $\bar
s\le s$, then $(p,q,s)_\wz$-atoms are $(p,\bar q,\bar
s)_\wz$-atoms, which further implies that
$h^{p,q,s}_{\wz}(\rz)\subset h^{p,\bar q,\bar s}_{\wz}(\rz)$ and
the inclusion is continuous.

Next we give some basic properties of $h^p_{\wz,N}(\rz)$ and
$h^{p,q,s}_{\wz}(\rz)$
\begin{prop}\label{p3.1.}\hspace{-0.1cm}{\rm\bf 3.1.}\quad
Let $\wz\in A_\fz^{loc}$. If $p\in (0,1]$ and $N\ge N_{p,\wz}$, then
the inclusion $h^p_{\wz,N}(\rz)\hookrightarrow {\cal S}'(\rz)$ is
continuous.
\end{prop}
Proof. Let $f\in h^p_{\wz,N}(\rz)$. For any $\vz\in {\cal
D}^0_N(\rz)$, and $\supp\vz\subset B_0=B(0,1)$, we have
$$\begin{array}{cl}
|<f,\vz>|&=|f*\bar\vz(0)|\le \|\bar \vz\|_{{\cal D}_N}\dinf_{x\in
B_0} {\cal M}^0_N(f)(x)\\
&\le [\wz(B_0)]^{-1/p}\|\vz\|_{{\cal
D}^0_N(\rz)}\|f\|_{h^p_{\wz,N}(\rz)},
\end{array}$$
where  $\bar\vz(x)=\vz(-x)$. This implies $f\in {\cal D}'(\rz)$
and the inclusion is continuous. The proof is finished.
\begin{prop}\label{p3.2.}\hspace{-0.1cm}{\rm\bf 3.2.}\quad
Let $\wz\in A_\fz^{loc}$. If $p\in (0,1]$ and $N\ge
[(n(q_\wz/p-1)]+2$, the the space $h^p_\wz(\rz)$ is complete.
\end{prop}
 Proof. For every $\vz\in {\cal D}_N^0(\rz)$ and every sequence
 $\{f_i\}_{i\in\nn}$ in ${\cal D}'(\rz)$ such that $\sum_i f_i$ converges
 in ${\cal D}'$ to the distribution $f$, the series $\sum_i
 f_i*\vz(x)$ converges pointwise to $f*\vz(x)$ for each $x\in\rz$.
 Thus,
$${\cal M}^0_Nf(x)^p\le \l(\dsum_i {\cal M}^0_N f_i(x)\r)^p\le
\dsum_i({\cal M}^0_Nf_i(x))^p\quad {\rm for \ all}\ x\in\rz,$$ and
hence $\|f\|_{h^p_{\wz,N}(\rz)}\le \dsum_i
\|f_i\|_{h^p_{\wz,N}(\rz)}$.

To prove that $h^p_{\wz,N}(\rz)$ is complete, it  suffices to show
that for every sequence $\{f_j\}_{j\in\nn}$ with
$\|f_j\|_{h^p_{\wz,N}(\rz)}<2^{-j}$ for any $j\in\nn$, the series
$\sum_{j\in\nn}f_j$ convergence in $h^p_{\wz,N}(\rz)$. Since
$\{\sum_{i=1}^j f_i\}_{j\in N}$ are Cauchy sequences in $h^p_{\wz,
N}(\rz)$, by Proposition 3.1 and the completeness of ${\cal
D}'(\rz)$, $\{\sum_{i=1}^j f_i\}_{j\in\nn}$ are also Cauchy
sequences in ${\cal D}'(\rz)$ and thus converge to some $f\in
{\cal D}'(\rz)$. Therefore,
$$\|f-\dsum_{i=1}^j f_i\|^p_{h^p_{\wz,N}(\rz)}=\|\dsum_{i=j+1}^\fz f_i\|^p_{h^p_{\wz,N}(\rz)}
\le \dsum_{i=j+1}^\fz 2^{-ip}\to 0$$ as $j\to\fz$. This finishes the
proof.

\begin{thm}\label{t3.1.}\hspace{-0.1cm}{\rm\bf 3.1.}\quad
Let $\wz\in A_\fz^{loc}$. If $(p,q,s)_\wz$ is an admissible
triplet and $N\ge N_{p,\wz}$, then $h^{p,q,s}_{\wz}(\rz)\subset
h^{p,q,s}_{\wz,N_{p,\wz}}(\rz)\subset h^p_{\wz,N}(\rz)$, and
moreover, there exists a positive constant $C$ such that for all
$f\in h^{p,q,s}_{\wz}(\rz)$,
$$\|f\|_{h^p_{\wz,N}(\rz)}\le
\|f\|_{h^p_{\wz,N_{p,\wz}}(\rz)}\le
C\|f\|_{h^{p,q,s}_{\wz}(\rz)}.$$
\end{thm}
Proof. Obviously, we only need to prove $h^{p,q,s}_{\wz}\subset
h^p_{\wz,N_{p,\wz}}(\rz)$ for  all $f\in h^{p,q,s}_{\wz}(\rz)$ ,
$\|f\|_{h^p_{\wz,N_{p,\wz}}(\rz)}\le \|f\|_{h^{p,q,s}_{\wz}(\rz)}$
. To this end, it suffice to prove that there exists a positive
constant $C$ such that
$$\|{\cal M}^0_{N_{p,\wz}}(a)\|_{L^p_\wz(\rz)}\le C
\quad{\rm for\ all}\ (p,q,s)_\wz-{atoms}\ a, \eqno(3.1)$$ and
$$\|{\cal M}^0_{N_{p,\wz}}(a)\|_{L^p_\wz(\rz)}\le C
\quad{\rm for\ a\ }\ (p,q)_\wz\ {\rm single\ atoms}\ a.
\eqno(3.2)$$ Since $q\in (q_\wz,\fz]$, so $\wz\in A_q^{loc}$. We
first prove (3.2). Let $a$ is a $(p,q)_\wz$ single atom. Using the
H\"older inequality, the $L^q_\wz(\rz)$-boundedness of ${\cal
M}_{N_{p,\wz}}^0$ and $\wz\in A_q^{loc}$ together with Proposition
2.2 (i), we have
$$\|{\cal M}^0_{N_{p,\wz}}(a)\|_{L^p_\wz(\rz)}\le
C\|a\|_{L^q_\wz(\rz)}^p[\wz(\rz)]^{1-p/q}\le C.$$ It remains to
prove (3.1). Let $a$  be  a $(p,q,s)_\wz$-atom supported in
$Q=Q(x_0,r)$. The first case is when $|Q|<1$. Then if $\bar Q$ is
the double of $Q$,
$$\begin{array}{cl}\dint_\rz [M_{N_p,\wz}^0(a)(x)]^p\wz(x)dx&=
\dint_{\bar Q} [M_{N_p,\wz}^0(a)(x)]^p\wz(x)dx+\dint_{\bar Q^c}
[M_{N_p,\wz}^0(a)(x)]^p\wz(x)dx\\
&:=I_1+I_2.
\end{array} $$
For $I_1$, by the properties of $A_q^{loc}$ (see Lemma 2.1), we have
$$I_1\le
C\|a\|_{L^q_\wz(\rz)}^p[\wz(\bar Q)]^{1-p/q}\le C.$$ To estimate
$I_2$, we claim that for  $x\in \bar Q^c$
$$M_{N_p,\wz}^0(a)(x)\le
C|x-x_0|^{s_0+1+n}|Q|^{s_0/n}[\wz(Q)]^{-1/p}\chi_{\{|x-x_0|<4n\}}(x),\eqno(3.3)$$
where $s_0=[n(q_\wz/p-1)]$. Indeed, let $P$ be the Taylor
expansion of $\vz$ at the point $(x-x_0)/t$ of order $s_0$. Thus,
by the Taylor remainder theorem, note that $0<t<1$, we then have
$$\begin{array}{cl}
|(a*\vz_t)(x)(x)|&=\l|t^{-n}\dint_\rz a(y)\l(\vz\l(\frac
{x-y}t\r)-P\l(\frac
{x_0-y}t\r)\r)dy\r|\\
&\le
C\chi_{\{|x-x_0|<4n\}}(x)|x-x_0|^{-(s_0+n+1)}\dint_B|a(y)||y|^{s_0+1}dy\\
&\le
C|x-x_0|^{s_0+1+n}|Q|^{(s_0+1)/n}[\wz(Q)]^{-1/p}\chi_{|x-x_0|<4n}(x).
\end{array} $$
Hence, (3.3) holds.  Choose $\eta>0$ such that, then by $\wz\in
A_{q_\wz+\eta}^{loc}$ and Proposition 2.2 (i), we have
$$I_2\le
C|Q|^{p(n+s_0/n)}[\wz(Q)]^{-1}\dint_{2r<|x-x_0|<4n}|x-x_0|^{p(s_0+1+n)}\wz(x)dx
\le C.$$ To deal with the case when $|Q|\ge 1$, the proof is
simple. In fact, let $Q^*=Q(x_0, r+n)$, by Corollary 2.1, we get
$$\begin{array}{cl}
\dint_\rz [M_{N_p,\wz}^0(a)(x)]^p\wz(x)dx&=\dint_{Q^*}
[M_{N_p,\wz}^0(a)(x)]^p\wz(x)dx\\
&\le C\|a\|_{L^q_\wz(\rz)}^p[\wz(Q^*)]^{1-p/q}\\
&\le C\|a\|_{L^q_\wz(\rz)}^p[\wz(Q)]^{1-p/q}\\
&\le C.\end{array}$$ Thus, Proposition 3.2 is proved.

\begin{center} {\bf 4. Calder\'on-Zygmund decompositions}\end{center}
In this section, we establish the Calder\'on-Zygmund
decompositions associated with grand maxiaml functions on weighted
$\rz$. We  follow the constructions in \cite {st}, \cite {b} and
\cite{blyz}.

Throughout  this section, we consider a distribution $f$ so that
for all $\lz>0$,
$$\wz(\{x\in\rz: {\cal M}_N(f)>\lz\})<\fz,$$ where $N\ge 2$ is some
fixed integer. Later with regard to the weighted local Hardy space
$h^p_{\wz,N}(\rz)$ with $p\in (0,1]$, we restrict to
$$N>[nq_\wz/p].$$
For a given $\lz>\dinf_{x\in\rz}{\cal M}_Nf(x)$, we set
$$\oz\equiv\{x\in\rz: {\cal M}_N(f)(x)>\lz\},$$
which  implies $\oz$ is a proper subset of $\rz$. As in
\cite{st1}, we give the usual Whitney decomposition of $\oz$. Thus
we can find closed cubes $Q_k$ whose interiors distance from
$\oz^c$, with $\oz=\cup_k Q_k$ and
$$diam(Q_k)\le 2^{-(6+n)} dist(Q_k,\oz)\le 4diam(Q_k).$$
Next, fix $a=1+2^{-(11+n)}$ and $b=1+2^{-(10+n)} $; if $\bar
Q_k=aQ_k, Q^*_k=bQ_k$, then $Q_k\subset\bar Q_k\subset Q_k^*$.
Also, $\bigcup Q_k^*=\oz$, and the $\{Q_k^*\}$ have the bounded
interior property: every point is contained in at most a fixed
number of the $\{Q_k^*\}$.

Fix a positive smooth function $\xi$ that equal $1$ in the cube of
side length $1$ centered at the origin and vanishes outside the
concentric cube of side length $a$. We set
$\xi_k(x)=\xi([x-x_k]/l_k)$, where $x_k$ is the center of the cube
$Q_k$ and $l_k$ is its side length. Obviously, for any $x\in\rz$,
we have $1\le \sum_k\xi_k(x)\le L$.  Write
$\eta_k=\xi_k/(\sum_j\xi_j)$. The $\eta_k$ form a partition of
unity for the set $\oz$ subordinate to the locally finite cover
$\{\bar Q_k\}$ of $Q$; that is to say, $\chi_\oz=\sum \eta_k$ with
each $\eta_k$ supported in the cube $Q_k$.

Let $s\in\nn_0$ be some fixed integers and ${\cal P}_s(\rz)$
denote the linear space of polynomials in $n$ variables of degrees
no more than $s$. For each $i$ and $P\in {\cal P}_s(\rz)$, set
$$\|P\|_i\equiv\l[\dfrac 1{\dint_\rz
\eta_i(x)dx}\dint_\rz|P(x)|^2\eta_i(x)dx\r]^{1/2}.\eqno(4.1)$$
Then $({\cal P}_s(\rz),\|\cdot\|_i)$ is a finite dimensional
Hilbert space . Let $f\in {\cal D}'(\rz)$. Since $f$ induces a
linear functional on ${\cal P}_s(\rz)$ via $Q|\to1/\int_\rz
\eta_i(x)dx<f,Q\eta_i>$, by the Riesz lemma, there exists a unique
polynomial $P_i\in {\cal P}_s(\rz)$ for each $i$ such that for all
$Q\in {\cal P}_s(\rz)$,
$$\begin{array}{cl}
\dfrac 1{\int_\rz \eta_i(x)dx}<f,Q\eta_i>&= \dfrac 1{\int_\rz
\eta_i(x)dx}<P_i,Q\eta_i>\\
&= \dfrac 1{\int_\rz \eta_i(x)dx}\dint_\rz P_i(x)Q(x)\eta_i(x)dx.
\end{array}$$
For every $i$, define distribution $b_i=(f-P_i)\eta_i$ if $l_i<1$,
 we set  $b_i=f\eta_i$ if $l_i\ge 1$.

We will show that for suitable choices of $s$ and $N$, the series
$\sum_i b_i$ converges in ${\cal D}'(\rz)$, and in this case, we
define $g=f-\sum_i b_i$ in ${\cal D}'(\rz)$.

The representation $f=g+\sum_i b_i$, where $g$ and $b_i$ are as
above, is said to be a Calder\'on-Zygmund decomposition of degree
$s$ and the height $\lz$ associated with ${\cal M}_N(f)$.

The rest of this section consists of series of lemmas. In Lemma
4.1 and Lemma 4.2, we give some properties of the smooth partition
of unity $\{\eta_i\}_i$. In Lemmas 4.3 through 4.6, we derive some
estimates for the bad parts $\{b_i\}_i$. Lemma 4.7 and Lemma 4.8
give controls over the good part $g$. Finally, Corollary 4.1 shows
the density of $L_\wz^q(\rz)\bigcap h^p_{\wz,N}(\rz)$ in
$h^p_{\wz,N}(\rz)$, where $q\in (q_\wz,\fz)$.
\begin{lem}\label{l4.1.}\hspace{-0.1cm}{\rm\bf 4.1.}\quad
There exists a positive constant $C_1$, depending only on $N$,
such that for all $i$ and $l\le l_i$,
$$\dsup_{|\az|\le N}\dsup_{x\in\rz}|\pz^\az\eta_i(lx)|\le
C_1.$$
\end{lem}
Lemma 4.1 is essentially Lemma 5.2 in \cite{b}.
\begin{lem}\label{l4.2.}\hspace{-0.1cm}{\rm\bf 4.2.}\quad
If  $l_i<1$, then there exists a constant a constant $C_2>0$
independent of $f\in {\cal D}'(\rz)$, $l_i$ and $\lz>0$ so that
$$\dsup_{y\in\rz}|P_i(y)\eta_i(y)|\le C_2\lz.$$
\end{lem}
Proof.  As in the proof of Lemma 5.3 in \cite{b}. Let
$\pi_l,\cdots, \pi_m(m=dim{\cal P}_s)$ be an orthonormal basis of
${\cal P}_s$ with respect to the norm (4.1). we have
$$P_i=\dsum_{k=1}^m\l(\dfrac 1{\int\eta_i}\dint
f(x)\pi_k(x)\eta_i(x)dx\r)\bar\pi_k,\eqno(4.2)$$ where the
integral is understood as $<f,\pi_k\eta_i>$. Hence
$$\begin{array}{cl}
1&=\dfrac 1{\int\eta_i}\dint|\pi_k(x)|^2\eta_i(x)dx\ge \dfrac
{2^{-n}}{|
Q_k|}\dint_{ Q_k}|\pi_k(x)|^2\eta_i(x)dx\\
&\ge  \dfrac {2^{-n}}{| Q_k|}\dint_{
Q_k}|\pi_k(x)|^2dx=2^{-n}\dint_{Q^0}|\wt\pi_k(x)|^2dx,
\end{array}\eqno(4.3)$$
where $\wt\pi_k(x)=\pi_k(x_i+l_i x)$ and $Q^0$ denotes the cube of
side length $1$ centered at the origin.

Since ${\cal P}_s$ is finite dimensional all norms on ${\cal P}_s$
are equivalent, there exists $A_1>0$ such that for all $p\in {\cal
P}_s$
$$\dsup_{|\az|\le s}\dsup_{z\in bQ^0}|\pz^\az P(z)|\le
A_1\l(\dint_{Q^0}|P(z)|^2dz\r)^{1/2}.$$ From this and (4.3), for
$k=1,\cdots, m$, we have
$$\dsup_{|\az|\le s}\dsup_{z\in bQ^0}|\pz^\az \wt\pi_k(z)|\le A_1.\eqno(4.4)$$
For $k=1,\cdots, m$ define
$$\Phi_k(y)=\dfrac{l_i}{\int\eta_i}\pi_k(z-l_iy)\eta_i(z-l_iy),$$
where $z$ is some point in $2^{9+n}nQ_k\bigcap\oz^c$.

It is easy to see that $\supp \Phi_k\subset B_n:=B(0,
2^{3(10+n)})$ and $\|\Phi_k\|_{{\cal D}_N}\le A_2$ by Lemma 4.1.

 Note that
$$\dfrac 1{\int\eta_i}\dint
f(x)\pi_k(x)\eta_i(x)dx=(f*(\Phi_k)_{l_i})(z),$$ since $l_i<1$, we
then have
$$\l|\dfrac 1{\int\eta_i}\dint
f(x)\pi_k(x)\eta_i(x)dx\r|\le {\cal M}_Nf(z)\|\Phi_k\|_{{\cal
D}_N}\le A_2\lz. $$ By (4.2), (4.4) and above estimate
$$\dsup_{z\in Q^*_i}| P_i(z)|\le
mA_1A_2\lz.$$ Thus, $$\dsup_{z\in\rz}| P_i(z)\eta_i(z)|\le
C_2\lz.$$ The proof is complete.
\begin{lem}\label{l4.3.}\hspace{-0.1cm}{\rm\bf 4.3.}\quad
There exists a constant $C_3>0$ such that
$${\cal M}_N^0b_i(x)\le C_3 {\cal M}_Nf(x)\quad {\rm for}\ x\in Q^*.\eqno(4.5)$$
\end{lem}
Proof. Take $\vz\in {\cal D}^0_N$,  and $x\in Q^*_i$.

Case I. For $t\le l_i$, we write
$$(b_i*\vz_t)(x)=(f*\Phi_t)(x)-((P_i\eta_i)*\vz_t)(x),$$
where $\Phi(z):=\vz(z)\eta_i(x-tz)$. Define
$\bar\eta_i(z)=\eta_i(x-l_iz))$. Obviously, $\supp \Phi\subset
B_n$. By Lemma 4.1, there exists a positive constant $C$ such that
$$\|\Phi\|_{{\cal D}_N}\le A_1C.$$
Note that for $N\ge 2$ there is a constant $C>0$ so that
$\|\vz\|_{L^1(\rz)}\le C$ for all $\vz\in {\cal D}^0_N$.
Therefore, by Lemma 4.2 and (4.5), we have
$$|b_i*\vz_l(x)|\le \|\Phi\|_{{\cal D}_N}{\cal M}_N
f(x)+C_2\lz\|\vz\|_{L^1(\rz)}\le C_3{\cal M}_N f(x),$$ since
${\cal M}_N f(x)>\lz$ for $x\in\oz$.

Case II. For $l_i<t<1$ by a simple calculation we can write
$$(b_i*\vz_t)(x)=\dfrac
{l_i}t(f*\Phi_{l_i})(x)-((P_i\eta_i)*\vz_t)(x),$$ where
$\Phi(z)=\vz(l_iz/t)\eta_i(x-l_i z)$. Define
$\bar\vz(z):=\vz(l_iz/t)$ and $\bar\eta_i(z)=\eta_i(x-l_iz)$. It
is easy to see that $\supp\Phi\subset B_n$. By Lemma 4.1, we can
find a  positive constant $C$ independent of $1>t>l_i$ so that
$$\dsup_{|\az|\le N}\dsup_{z\in\rz}|\pz^\az\bar\vz(z)|\le C,\quad
\dsup_{|\az|\le N}\dsup_{z\in\rz}|\pz^\az\bar\eta_i(z)|\le A_1.$$
Hence,  there exists  a positive constant $C$ such that
$\|\Phi\|_{{\cal D}_N} \le C$, and $\|\vz\|_{L^1(\rz)}\le C$ for
$\vz\in{\cal D}_N$ for $N\ge 2$. As in the case I
$$|(b_i*\vz_t)(x)|\le \|\Phi\|_{{\cal D}_N} {\cal
M}_Nf(x)+C_2\lz\|\vz\|_{L^1(\rz)}\le C{\cal M}_Nf(x).$$ By
combining both cases, we can obtain the desired result.
\begin{lem}\label{l4.4.}\hspace{-0.1cm}{\rm\bf 4.4.}\quad
Suppose $Q\subset \rz$ is bounded, convex, and $0\in Q$, and $N$
is a positive integer. Then there is a constant $C$ depending only
on $Q$ and $N$ such that for every $\phi\in{\cal D}(\rz)$ and
every integer $s,\ 0\le s<N$ we have
$$\dsup_{x\in Q}\dsup_{|\az|\le N}|\pz^\az R_y(z)|\le C
\dsup_{x\in Q}\dsup_{s+1\le |\az|\le N}|\pz^\az \phi(z)|,$$ where
$R_y$ is the remainder of the Taylor expansion of $\phi$ of order
$s$ at the point $y\in\rz$.
\end{lem}
Lemma 4.4 is Lemma 5.5 in \cite{b}.

\begin{lem}\label{l4.5.}\hspace{-0.1cm}{\rm\bf 4.5.}\quad
Suppose $0\le s<N$. Then there exist positive constants $C_3, C_4$
so that for $i\in\nn$,
$${\cal M}_N^0(b_i)(x)\le
C\dfrac{l_i^{n+s+1}}{(l_i+|x-x_i|)^{n+s+1}}\chi_{\{|x-x_i|<C_3\}}(x)\
{\rm if}\ x\not\in Q_i^*.\eqno(4.6)$$ Moreover,
$${\cal M}_N^0(b_i)(x)=0,\ {\rm if}\  x\not\in Q_i^*\ {\rm and}\
l_i\ge C_4.$$
\end{lem}
Proof. Take $\vz\in {\cal D}(\rz)$. Recall that $\eta_i$ is
supported in the cube $\bar Q_i$, and we have taken $\bar Q_i$ to
be strictly contained in $Q_i^*$. Thus if $x\not\in Q_i^*$ and
$\eta_i(y)\not=0$, then there exists a positive constant $C_3$
such that  $|x-y|\le|x-x_i|\le C_3|x-y|$, and the support property
of $\Phi$ requires that $1>t\ge |x-y|\le 2^{-11-n}l_i$. Hence,
$|x-x_i|\le C_3t$ and $l_i< 2^{11+n}:=C_4$ and $l_i<C_4t$ . Pick
some $w\in (2^{8+n}n Q_i)\bigcap\oz^c$.

 Case I. If $1\le l_i< C_4$ and $\vz\in {\cal D}^0_N$, where  define
 $\phi(z)=\vz(\bar l_iz/t)$ and $\bar l_i=l_i/C_4< 1$.  We have
 $$\begin{array}{cl}(b*\vz_l)(x)&=t^{-n}\dint b_i\vz((x-z)/t)dz\\
 &=t^{-n}\dint b_i\phi((x-z)/\bar l_i)dz\\
&=t^{-n}\dint b_i\phi_{(x-w)/\bar l_i}((w-z)/\bar l_i)dz\\
&=\dfrac{\bar l_i^n}{t^n}(f*\Phi_{\bar l_i})(w), \end{array}$$
where
$$\Phi(z):=\phi_{(x-w)/\bar l_i}(z)\eta_i(w-\bar l_iz), \quad \phi_{(x-w)/\bar l_i}(z)=\phi(z+(x-w)/\bar l_i).$$
Obviously, $\supp\Phi\subset B_n$. Note that $l_i<tC_4$ and
$|x-x_i|\le C_3t$, we obtain
$$|(b*\vz_t)(x)|\le C\dfrac{\bar l_i^n}{t^n}{\cal M}_Nf(w)\le C\lz \dfrac{\bar l_i^n}{t^n}\le
C\lz\dfrac{l_i^{n+s+1}}{(l_i+|x-x_i|)^{n+s+1}}.\eqno(4.7)$$
 Case II. If $ l_i<1$ and $\vz\in {\cal D}^0_N$ define
 $\phi(z)=\vz(l_iz/t)$.  Consider the Taylor expansion of $\phi$ of order $s$
at the point $y:=(x-w)/l_i$,
$$\phi(y+z)=\dsum_{|\az|\le
s}\dfrac{\pz^\az\phi(y)}{\az!}z^\az+R_y(z),$$ where $R_y$ denotes
the remainder.

Thus,
 $$\begin{array}{cl}(b*\vz_t)(x)&=t^{-n}\dint b_i\vz((x-z)/t)dz\\
 &=t^{-n}\dint b_i\phi((x-z)/l_i)dz\\
&=t^{-n}\dint b_iR_{(x-w)/l_i}((w-z)/l_i)dz\\
&=\dfrac{l_i^n}{t^n}(f*\Phi_{l_i})(w)+t^{-n}\dint
P_i(z)\eta_i(z)R_{(x-w)/l_i}((w-z)/l_i)dz,
\end{array}\eqno(4.8)$$ where
$$\Phi(z):=R_{(x-w)/l_i}(z)\eta_i(w-l_iz).$$
Obviously, $\supp\Phi\subset B_n$. Apply Lemma 4.4 to
$\phi(z)=\vz(l_iz/t),\ y=(x-w)/l_i$ and $Q=B_n$. We have
$$\begin{array}{cl}
\dsup_{z\in B_n}\dsup_{|\az|\le N}|\pz^\az R_y(z)|&\le C
\dsup_{z\in y+B_n}\dsup_{s+1\le |\az|\le N}|\pz^\az \phi(z)|\\
&\le C\dsup_{z\in y+B_n}\l(\dfrac {l_i}t\r)^{-(s+1)}\dsup_{s+1\le |\az|\le N}|\pz^\az \vz(l_iz/t)|\\
&\le C\l(\dfrac {l_i}t\r)^{-(s+1)}.
\end{array}$$
 Note that $l_i<tC_4$ and
$|x-x_i|\le C_3t$,  therefore by (4.8), we have
$$\begin{array}{cl}
(b*\vz_t)(x)&\le\dfrac{l_i^n}{t^n}|(f*\Phi_{l_i})(w)|+t^{-n}\dint
|P_i(z)\eta_i(z)R_{(x-w)/l_i}((w-z)/l_i)|dz\\
&\le C\l(\dfrac{l_i^n}{t^n}{\cal M}_Nf(w)\|\Phi\|_{{\cal D}_N}+
\lz\dsup_{z\in
B_n}\dsup_{|\az|\le N}|\pz^\az R_y(z)|\r]\\
&\le C\lz\dfrac{l_i^{n+s+1}}{(l_i+|x-x_i|)^{n+s+1}}.
\end{array}\eqno(4.9)$$
Combining (4.7) and (4.9), we obtain (4.6).

\begin{lem}\label{l4.6.}\hspace{-0.1cm}{\rm\bf 4.6.}\quad
Let $\wz\in A_\fz^{loc}$ and $q_\wz$ be as in (2.3). If
$p\in(0,1]$, $s\ge [nq_\wz/p]$ and $N>s$, there exists a positive
constant $C_5$ such that for all $f\in h^{p}_{\wz,N}(\rz)$,
$\lz>\dinf_{x\in\rz}{\cal M}_Nf(x)$ and $i$,
$$\dint_\rz[{\cal M}^0_N(b_i)(x)]^p\wz(x)dx\le C_5\dint_{Q_i^*}
[{\cal M}_N(f)(x)]^p\wz(x)dx.\eqno(4.10)$$ Moreover the series
$\sum_i b_i$ converges in $h^p_{\wz,N}(\rz)$ and
$$\dint_\rz[{\cal M}^0_N(\sum_i b_i)(x)]^p\wz(x)dx\le  C_5\dint_{\oz}
[{\cal M}_N(f)(x)]^p\wz(x)dx.\eqno(4.11)$$
\end{lem}
Proof. By Lemma 4.4, we have
$$\begin{array}{cl}
\dint_\rz[{\cal M}^0_N(b_i)(x)]^p\wz(x)dx&\le \dint_{Q^*_i}[{\cal
M}^0_N(b_i)(x)]^p\wz(x)dx\\
&\qquad+\dint_{C_3Q^0_i\setminus Q^*_i}[{\cal
M}^0_N(b_i)(x)]^p\wz(x)dx,
\end{array}\eqno(4.12)$$
where $Q_i^0=Q(x_i,1)$. Note that  $s\ge[nq_\wz/p]$ implies
$2^{-n(q_\wz+\eta)}2^{(s+n+1)p}>1$ for sufficient small $\eta>0$.
Using Lemma 2.1 (ii) with $\wz\in A_{q_\wz+\eta}^{\loc}$, Lemma
4.5 and  the fact that ${\cal M}_N(f)(x)>\lz$ for all $x\in
Q^*_i$, we have
$$\begin{array}{cl}
\dint_{C_3Q^0_i\setminus Q^*_i}[{\cal
M}^0_N(b_i)(x)]^p\wz(x)dx&\le
\dsum_{k=0}^{k_0}\dint_{2^kQ_i^*\setminus 2^{k-1}Q_i^*}[{\cal
M}^0_N(b_i)(x)]^p\wz(x)dx\\
&\le
\lz^p\wz(Q^*_i)\dsum_{k=0}^{k_0}[2^{-n(q_\wz+\eta)+(s+n+1)p}]^{-k}\\
&\le C\dint_{Q_i^*}[{\cal M}_Nf(x)]^p\wz(x)dx,
\end{array}\eqno(4.13)$$
where $k_0\in\zz$ such that $2^{k_0-1}\le C_3< 2^{k_0}$.

Combining (4.12) and (4.13), then (4.10) holds.
 By (4.10), we have
 $$\dint_\rz[{\cal M}^0_N(b_i)(x)]^p\wz(x)dx\le C
 \dsum_i\dint_{Q_i^*}[{\cal M}_Nf(x)]^p\wz(x)dx\le C\dint_\oz
[{\cal M}_N(f)(x)]^p\wz(x)dx,$$ which together with complete of
$h^p_{\wz,N}$ (see Proposition 3.2)implies that $\sum_i b_i$
converges in $h^p_{\wz,N}$. So by Proposition 3.1, the series
$\sum_i b_i$ converges in ${\cal D}'(\rz)$, and therefore ${\cal
M}_N^0(\sum_i b_i)(x)\le \sum_i {\cal M}_N^0( b_i)(x)$, which
gives (4.11). Thus, Lemma 4.6 is proved.

\begin{lem}\label{l4.7.}\hspace{-0.1cm}{\rm\bf 4.7.}\quad
Let $\wz\in A_\fz^{loc}$ and $q_\wz$ be as in (2.3), $s\in\nn_0$,a
nd $N\ge 2$. If $q\in(q_\wz,\fz]$ and $f\in L_\wz^q(\rz)$, then
the series $\sum_i b_i$ converges in $L_\wz^q(\rz)$ and there
exists a positive constant $C_6$, independent of $f$ and $\lz$,
such that $\|\sum_i |b_i|\|_{L^q_\wz(\rz)}\le
C_6\|f\|_{L^q_\wz(\rz)}$.
\end{lem}
Proof. The proof for $q=\fz$ is similar to that for $q\in
(q_\wz,\fz)$. So we only give the proof for $q\in(q_\wz,\fz)$. Set
$F_1=\{i\in\nn: |Q_i|\ge 1\}$ and $F_2=\{i\in\nn: |Q_i|< 1\}$. By
lemma 4.3, for $i\in F_2$, we have
$$\begin{array}{cl}
\dint_\rz|b_i(x)\wz(x)dx&\le\dint_{Q_i^*}|f(x)|^q\wz(x)dx+
\dint_{Q_i^*}|P_i(x)\eta_i(x)|^q\wz(x)dx\\
&\le \dint_{Q_i^*}|f(x)|^q\wz(x)dx+ \lz^q\wz(Q^*_i).\end{array}$$
For $i\in F_1$, we have
$$\dint_\rz|b_i(x)|\wz(x)dx\le\dint_{Q_i^*}|f(x)|^q\wz(x)dx.$$
From these, we obtain
$$\begin{array}{cl}
\dsum_i\dint_\rz|b_i(x)|\wz(x)dx&=\dsum_{i\in
F_1}\dint_\rz|b_i(x)|\wz(x)dx+\dsum_{i\in
F_2}\dint_\rz|b_i(x)|\wz(x)dx\\&\le\dint_{Q_i^*}|f(x)|^q\wz(x)dx+
\dint_{Q_i^*}|P_i(x)\eta_i(x)|^q\wz(x)dx\\
&\le \dsum_i\dint_{Q_i^*}|f(x)|^q\wz(x)dx+ C\dsum_{i\in
F_2}\lz^q\wz(Q^*_i)\\
&\le \dsum_i\dint_{Q_i^*}|f(x)|^q\wz(x)dx+ C\lz^q\wz(Q)\\
&\le C_6\dint_\rz|f(x)|^q\wz(x)dx.
\end{array}$$
From this and applying $b_i$ have finite covers, we have
$$\|\sum_i |b_i|\|_{L^q_\wz(\rz)}\le C_6\|f\|_{L^q_\wz(\rz)}.$$
The proof is finished.

\begin{lem}\label{l4.8.}\hspace{-0.1cm}{\rm\bf 4.8.}\quad
 If $N>s\ge 0$ and $\sum_i b_i$ converges in ${\cal D}'(\rz)$,
 then there exists a  positive constant $C_7$, independent of $f$
 and $\lz$, such that for all $x\in\rz$,
 $${\cal M}_N^0(g)(x)\le {\cal
M}_N^0(f)(x)\chi_{\oz^c}(x)+
C_7\dfrac{l_i^{n+s+1}}{(l_i+|x-x_i|)^{n+s+1}}\chi_{\{|x-x_i|<C_3\}}(x)$$
\end{lem} Proof. If $x\not\in\oz$, since ${\cal M}_N^0(g)(x)\le {\cal
M}_N^0(f)(x)+\dsum_i {\cal M}_N^0(b_i)(x)$,  by Lemma 4.5, we
obtain
$${\cal M}_N^0(g)(x)\le {\cal
M}_N^0(f)(x)\chi_{\oz^c}(x)+
C\dsum_i\dfrac{l_i^{n+s+1}}{(l_i+|x-x_i|)^{n+s+1}}\chi_{\{|x-x_i|<C_3\}}(x).$$
If $x\in\oz$, choose $k\in\nn$ such that $x\in Q_k^*$. Let
$J:=\{i\in\nn: Q_i^*\bigcap Q_k^*\not=\O\}$. Then the cardinality
of $J$ is bounded by $L$. By Lemma 4.5, we have
$$\dsum_{i\not\in J} {\cal M}_N^0(b_i)(x)\le C\lz \dsum_{i\not\in J}
\dfrac{l_i^{n+s+1}}{(l_i+|x-x_i|)^{n+s+1}}\chi_{\{|x-x_i|<C_3\}}(x).$$
It suffices to estimate the grand maximal function of
$g+\sum_{i\not=J}b_i=f-\sum_{i\in J}b_i$.
  Take $\vz\in {\cal D}^0_N$ and $0<t<1$. We write
$$\begin{array}{cl}
(f-\sum_{i\in J}b_i)*\vz_t(x)&=(f\xi)*\vz_t+(\sum_{i\in
J}P_i\eta_i)*\vz_t\\
&=f*\Phi_t(w)+(\sum_{i\in J}P_i\eta_i)*\vz_t,
\end{array}$$
where $w\in (2^{8+n}nQ_k)\bigcap\oz^c$, $\xi=1-\sum_{i\in
J}\eta_i$ and
$$\Phi(z):=\vz(z+(x-w)/t)\xi(w-tz).$$
 Since for $N\ge 2$ there is a constant $C>0$ so that $\|\vz\|_{L^1(\rz)}\le
 C$ for all $\vz\in {\cal D}^0_N$ and Lemma 4.1, we have
 $$\l|\l(\dsum_{i\in J}P_i\eta_i)\r)*\vz_t(x)\r|\le C\lz.$$
 Finally, we estimate $f*\Phi_t(w)$. There are two cases: If $t\le
 2^{-(11+n)}l_k$,  then $f*\Phi_t(w)=0$, because $\xi$ vanishes in $Q_k^*$
 and $\vz_t$ is supported in $B(0,t)$. On the other hand, if $t\ge
 2^{-(11+n)}l_k$, then there exists a positive constant $C$ such that
 $\supp\Phi\subset B_n$ and $\|\Phi\|_{{\cal D}_N}\le
 C$. Hence,
 $$|(f*\Phi_t)|\le {\cal M}_N f(w)\|\Phi\|_{{\cal D}_N}\le C\lz.$$
 By the above estimates, we have
 $$|(f-\sum_{i\in J}b_i)*\vz_t|\le C\lz.$$
That is,
 $${\cal M}^0_N((f-\sum_{i\in J}b_i))(x)\le C\lz.$$
Thus, Lemma 4.8 is proved.
\begin{lem}\label{l4.9.}\hspace{-0.1cm}{\rm\bf 4.9.}\quad
Let $\wz\in A_\fz^{loc}$, $q_\wz$ be as in (2.3) and $p\in (0,1]$.
\begin{enumerate}
 \item[(i)] If $N>s\ge [nq_\wz/p]$ and ${\cal M}_N(f)\in
 L_\wz^p(\rz)$, then ${\cal M}_N(g)\in  L_\wz^1(\rz)$ and there
 exists a positive constant $C_8$, independent of $f$ and $\lz$,
 such that
 $$\dint_\rz[{\cal M}_N^0(g)(x)]^q\wz(x)dx\le C_8\lz^{1-p}\dint_\rz
[{\cal M}_N(f)(x)]^p\wz(x)dx.$$
 \item[(ii)] If $N\ge 2$ and $f\in L_\wz^1(\rz)$, then $g\in
 L_\wz^\fz(\rz)$ and there exists a positive constant $C_9$,
independent of $f$ and $\lz$,
 such that $\|g\|_{L_\wz^\fz}\le C_9\lz$.
\end{enumerate}
\end{lem}

Proof. Since $f\in h^p_{\wz,N}(\rz)$, by Lemma 4.6, $\sum_i b_i$
converges in $h^p_{\wz, N}(\rz)$ and there in ${\cal D}'(\rz)$ by
proposition 3.1. Observe that $s\ge[nq_\wz/p]$, by Lemma 4.8, we
obtain
$$\begin{array}{cl}
\dint_\rz [{\cal M}^0_N(g)(x)]\wz(x)dx&\le  C\lz\dsum_i\dint_\rz
\dfrac{l_i^{(n+s+1)}}{(l_i+|x-x_i|)^{(n+s+1)}}\chi_{\{|x-x_i|<C_3\}}(x)\wz(x)dx\\
&\qquad+\dint_{\oz^c}[{\cal M}_N(f)(x)]\wz(x)dx\\
&\le C\lz^q\dsum_i\wz(Q_i^*)+\dint_{\oz^c}[{\cal M}_N(f)(x)]\wz(x)dx\\
&\le C\lz\wz(\oz)+C\lz^{1-p}\dint_{\oz^c}[{\cal M}_N(f)(x)]^p\wz(x)dx\\
&\le C\lz^{1-p}\dint_{\oz^c}[{\cal M}_N(f)(x)]^p\wz(x)dx.
\end{array}$$
Thus, (i) holds.

Moreover, if $f\in L_\wz^1(\rz)$, then $g$ and $\{b_i\}$ are
functions, and Lemma 4.7, $\sum_i b_i$ converges in $L^q_\wz(\rz)$
and thus in ${\cal D}'(\rz)$ by Lemma 2.4. Write
$$g=f-\dsum_i b_i=f(1-\dsum_i\eta_i)+\dsum_{i\in
F_2}P_i\eta_i=f\chi_{\oz^c}+\dsum_{i\in F_2}P_i\eta_i.$$ By Lemma
4.3, we have $|g(x)|\le C\lz$ for all $x\in\oz$, and by
Proposition 2.3, $|g(x)|=|f(x)|\le {\cal M}_N f(x)\le \lz$ for
almost everywhere $x\in \oz^c$, which leads to that
$\|g\|_{L_\wz^\fz(\rz)}\le C\lz$ and thus yields (ii). The proof
is finished.

\begin{cor}\label{c4.1.}\hspace{-0.1cm}{\rm\bf 4.1.}\quad
Let $\wz\in A_\fz^{loc}$ and $q_\wz$ be as in (2.3). If $q\in
(q_\wz,\fz)$, $N>[nq_\wz/p]$ and $p\in (0,1]$, then
$h^p_{\wz,N}(\rz)\bigcap L_\wz^1(\rz)$ is dense in
$h^p_{\wz,N}(\rz)$.
\end{cor}
Proof. Let $f\in h^p_{\wz,N}(\rz)$. For any
$\lz>\dinf_{x\in\rz}{\cal M}_Nf(x)$, let $f=g^\lz+\sum_i b_i^\lz$
be the Calder\'on-Zygmund decomposition of $f$ of degree $s$ with
$[nq_\wz/p]\le s<N$ and height $\lz$ associated to ${\cal M}_N f$.
By Lemma 4.6,
$$\|\dsum_i b_i^\lz\|_{h^p_{\wz,N}(\rz)}\le C\dint_{\{x\in\rz: {\cal
M}_Nf(x)>\lz\}}[{\cal M}_Nf(x)]^p\wz(x)dx.$$ Therefore, $g^\lz\to
f$ in $h^p_{\wz, N}(\rz)$ as $\lz\to\fz$. But by Lemma 4.9, ${\cal
M}_N (g^\lz)\in L_\wz^1(\rz)$, so by Proposition 2.2, $g^\lz\in
L_\wz^1(\rz)$.  Thus, Corollary 4.1 is proved.

\begin{center} {\bf 5. Weighted atomic decompositions of $h^p_{\wz,N}(\rz)$}\end{center}
We will follow the proof of atomic decomposition as presented by
Stein in \cite{st}.

In this section, we take $k_0\in\zz$ such that
$2^{k_0-1}\le\dinf_{x\in\rz}{\cal M}_Nf(x)<2^{k_0}$, if
$\dinf_{x\in\rz}{\cal M}_Nf(x)=0$, write $k_0=-\fz$. Let $\wz\in
A_\fz^{loc}, q_\wz$ be as in  (2.3), $p\in (0,1]$ and $N>s\equiv
[nq_\wz/p]$.
Let $f\in h^p_{\wz,N}(\rz)$. For each integer $k\ge
k_0$ consider the Calder\'on-Zygmund decomposition of $f$ of
degree $s$ and height $\lz=2^k$ associated to ${\cal M}_Nf$,
$$f=g^k+\dsum_{i\in\nn} b_i^k,$$ where
$$\oz^k:=\{x\in\rz: {\cal M}_Nf(x)>2^k\}, \ Q_i^k:=Q_{l_i^k}$$ and
$b_i^k:=(f-P_i^k)\eta_i^k$ if $l_i^k<1$ and  $b_i^k:=f\eta_i^k$ if $l_i^k\ge1$.

Recall that for
fixed $k\ge k_0$, $(x_i=x_i^k)_{i\in\nn}$ is a sequence in $\oz^k$
and $(l_i=l_i^k)_{i\in\nn}$  for $\oz=\oz^k$, $\eta_i=\eta_i^k$
given in Section 4 and $P_i=P_i^k$ is the projection of $f$ onto
${\cal P}_s$ with respect to the norm given in Section 4.

Define a polynomial $P_{ij}^{k+1}$ as an orthogonal projection of
$(f-P_j^{k+1})\eta_i^j$ on ${\cal P}_s$ with respect to the norm
$$\|P\|^2=\dfrac 1{\int_\rz
\eta_j^{k+1}}\dint_\rz|P(x)|^2\eta_j^{k+1}(x)dx,$$ that is
$P_{ij}^{k+1}$ is the unique element of ${\cal P}_s$ such that
$$\dint_\rz
(f(x)-P_j^{k+1}(x))\eta_i^k(x)Q(x)\eta_j^{k+1}(x)dx=\dint_\rz
P_{ij}^{k+1}(x)Q(x)\eta_j^{k+1}(x)dx.$$  For convenience we denote
$Q_i^{k*}=(1+2^{-(9+n)})Q_i^k$, $E_1=\{i\in\nn: |Q_i|\ge
1/(2^4n)\}$ and $E_2=\{i\in\nn: |Q_i|< 1/(2^4n)\}$,
$F_1=\{i\in\nn: |Q_i|\ge 1\}$ and $F_2=\{i\in\nn: |Q_i|< 1\}$.

There are two things we need to know about the polynomials
$P_{ij}^{k+1}$. First, $P_{ij}^{k+1}\not=0$ only if $
Q_i^{k*}\bigcap Q_j^{k+1*}\not=\o$; this follows directly from the
the definition of $P_{ij}^{k+1}$( since it involves
$\eta_i^{k+1}$, which is supported in $Q_i^{j+1*}$). More
precisely, we have the following results.
\begin{lem}\label{l5.1.}\hspace{-0.1cm}{\rm\bf 5.1.}\quad
Note that $\oz^{k+1}\subset \oz^k$, then
\begin{enumerate}
\item[(i)] If $Q_i^{k*}\bigcap Q_j^{k+1*}\not=\o$, then
$l_j^{k+1}\le 2^4\sqrt{n}l_i^k$ and $Q_j^{k+1*}\subset
2^6nQ_j^{k*}\subset \oz^k$. \item[(ii)]There exists a positive $L$
such that for each $j\in\nn$ the cardinality of $\{i\in\nn:\
Q_i^{k*}\bigcap Q_j^{k+1*}\not=\o$ is bounded by $L$.
\end{enumerate}
\end{lem}
  \begin{lem}\label{l5.2.}\hspace{-0.1cm}{\rm\bf 5.2.}\quad
If $ l_j^{k+1}<1$,
$$\dsup_{y\in\rz}|P_{ij}^{k+1}(y)\eta_j^{k+1}(y)|\le C2^{k+1}.\eqno(5.1)$$
\end{lem}

\begin{lem}\label{l5.3.}\hspace{-0.1cm}{\rm\bf 5.3.}\quad
 For every $k\in\zz$, $\dsum_{i\in\nn}(\dsum_{j\in F_2}
 P_{ij}^{k+1}\eta_j^{k+1})=0$, where the series converges pointwise
 and in ${\cal D}'(\rz)$.
\end{lem}
 Lemmas 5.1-5.3  can be proved by the methods in Lemmas 6.1-6.3 in
\cite{b}.

The following lemma establishes the weighted atomic decompositions
for a dense subspace of $h^p_{\wz, N}(\rz)$.
\begin{lem}\label{l5.4.}\hspace{-0.1cm}{\rm\bf 5.4.}\quad
Let $\wz\in A_\fz^{loc}$ and $q_\wz$ be as in (2.3). If
$ p\in(0,1]$, $s\ge [nq_\wz/p]$ and $N>s$, then
for any $f\in (L^1_\wz(\rz)\bigcap h^p_{\wz, N}(\rz))$, there
exists numbers $\lz_0$ and $\{\lz_i^k\}_{k\in\zz,i}\subset \cc$,
$(p,\fz,s)_\wz$-atoms $\{a_i^k\}_{k\in\zz,i}$ and single atom
$a_0$ such that
$$f=\dsum_{k\in\zz}\dsum_{i}\lz_i^ka_i^k+\lz_0a_0,$$ where the series converges
almost everywhere and in ${\cal D}'(\rz)$, moreover, there exists
a positive $C$, independent of $f$, such that
$\dsum_{k\in\zz,i}|\lz_i^k|^p+|\lz_0|^p\le
C\|f\|_{h^p_{\wz,N}(\rz)}$.
\end{lem}
Proof. Let $f\in (L^1_\wz(\rz)\bigcap h^p_{\wz, N}(\rz))$. We
first consider the case $k_0=-\fz$. For each $k\in\zz$, $f$ has a
Calder\'on-Zygmund decomposition of degree $s\ge [nq_\wz/p]$ and
height $2^k$ associated to ${\cal M}_N(f), f=g^k+\sum_i b_i^k$ as
above.  By Corollary 4.1 and Proposition 3.1, $g^k\to f$ in both
$h^p_{\wz,N}(\rz)$ and ${\cal D}'(\rz)$ as $k\to\fz$. By Lemma
4.9(ii), $\|g^k\|_{L_\wz^p(\rz)}\to 0$ as $k\to-\fz$, and
moreover, by Lemma 2.2 (ii), $g^k\to 0$ in ${\cal D}'(\rz)$ as
$k\to-\fz$. Therefore,
$$f=\dsum_{k=-\fz}^\fz(g^{k+1}-g^k) \eqno(5.2)$$ in ${\cal D}'(\rz)$.
Moreover, since $supp(\sum_i b_i^k)\subset\oz_k$ and
$\wz(\oz_k)\to 0$ as $k\to \fz$, then $g^k\to f$ almost everywhere
as $k\to\fz$. Thus, (5.2) also holds almost everywhere.

By Lemma 5.1 and
$\sum_i\eta_i^kb_j^{k+1}=\chi_{\oz_k}b_j^{k+1}=b_j^{k+1}$ for all
$j$, then
$\sum_i\eta_i^kb_j^{k+1}=\chi_{\oz_k}b_j^{k+1}=b_j^{k+1}$ for all
$j$,
$$\begin{array}{cl}
g^{k+1}-g^k&=\l(f-\dsum_j b_j^{k+1}\r)-\l(f-\dsum_i b_i^{k}\r)\\
&=\dsum_i b_i^k-\dsum_j b_j^{k+1}\\
&=\dsum_i\l[ b_i^k-\dsum_{j\in F^k_1} b_j^{k+1}\eta_i^k+\dsum_{j\in
F^k_2}b_j^{k+1}\eta_i^{k}\r]\\
&\equiv\dsum_i h_i^k.
\end{array}$$
where $F^k_1=\{i\in\nn: |Q^k_i|\ge 1\}$ and $F^k_2=\{i\in\nn: |Q^k_i|< 1\}$ and  the series converges in ${\cal D}'(\rz)$ and almost
everywhere. Furthermore, we rewrite $h_i^k$ into
$$h_i^k=f\chi_{(\oz_{k+1})^c}\eta_i^k-P_i^k\eta_i^k+\dsum_j
P_j^{k+1}\eta_i^k\eta_j^{k+1}+\dsum_{j\in F_2}
P_j^{k+1}\eta_j^{k+1}.$$ By proposition 2.2, $|f(x)|\le {\cal M}_N
f(x)\le 2^{k+1}$ for almost everywhere $x\in (\oz_{k+1})^c$, and
by Lemma 4.2 and (5.1),
$$\|h_i^k\|_{L^\fz_\wz(\rz)}\le C2^k\ {\rm for }\  i\in\nn.\eqno(5.3)$$
Next we consider three cases about $i$.

Case I. When $i\in F_1$, we have
$$h_i^k=f\eta_i^k+\dsum_{j\in F_1}f\eta_j^{k+1}\eta_i^k+
\dsum_{j\in F_2}(f-P_j^{k+1})\eta_j^{k+1}\eta_i^k+\dsum_{j\in
F_2}P_{ij}^{k+1}\eta_i^{k+1}.$$
Case II. When $i\in E_1\bigcap F_2$,
we have
$$h_i^k=(f-P_i^k)\eta_i^k+\dsum_{j\in F_1}f\eta_j^{k+1}\eta_i^k+
\dsum_{j\in F_2}(f-P_j^{k+1})\eta_j^{k+1}\eta_i^k+\dsum_{j\in
F_2}P_{ij}^{k+1}\eta_i^{k+1}.$$ Case III. When $i\in E_2$, if
$j\in F_1$, then $l_i^k<  l_j^{k+1}/(2^4n)$, so $ Q_i^{j*}\bigcap
Q_j^{k+1*}=\o$ by Lemma 5.1 (i). Thus,  we have
$$\begin{array}{cl}
h_i^k&=(f-P_i^k)\eta_i^k+\dsum_{j\in F_1}f\eta_j^{k+1}\eta_i^k+
\dsum_{j\in F_2}(f-P_j^{k+1})\eta_j^{k+1}\eta_i^k+\dsum_{j\in
F_2}P_{ij}^{k+1}\eta_i^{k+1}\\
&=(f-P_i^k)\eta_i^k+ \dsum_{j\in
F_2}(f-P_j^{k+1})\eta_j^{k+1}\eta_i^k+\dsum_{j\in
F_2}P_{ij}^{k+1}\eta_i^{k+1},
\end{array}$$
We next let $\gz=1+2^{-12-n}$.

For Cases I and II. Obviously, $h_i^k$ is supported in a cube $\wt
Q_i^k$ that contains $Q_i^{k*}$ as well as all the $Q_j^{k+1*}$
that intersect $Q_i^{k*}$. In fact, observe that if
$Q_i^{k*}\bigcap Q_j^{k+1*}\not=\o$, by Lemma 5.1, we have
 $$Q_j^{k+1*}\subset
2^6nQ_j^{k*}\subset \oz^k.$$ So, if $l^k_i<Ln/(\gz-1)$, we set
$$\wt Q_i^k:=2^6nQ_j^{k*}.$$
 On the other hand, note that $l_j^{k+1}<1$ and $ l_i^k\ge
2^{-n-4}$, then $Q_j^{k+1*}\subset Q(x^k_i, l^k_i+Ln)$. So, if
$l^k_i\ge Ln/(\gz-1)$, we set $\wt Q_i^k=\gz Q_j^{k}.$ Hence,
$$Q_j^{k+1*}\subset Q(x^k_i, l^k_i+Ln)\subset \wt Q_i^k=\gz Q_i^k= Q_i^{k*}\subset \oz^k,$$
 if $l^k_i\ge Ln/(\gz-1)$.

  From these, for Cases I and II, there exists  a positive constant $C_{10}$ such that
  $$ \wt Q_i^k\subset \oz^k,\ {\rm and}\ \wz(\wt Q_i^k)\le C_{10}\wz(Q_i^{k*}).$$
   But, $h_i^k$ does not satisfy the moment conditions.

For Case  III. We claim that $h_i^k$ is supported in a cube $\wt
Q_i^k$ that contains $Q_i^{k*}$ as well as all the $Q_j^{k+1*}$
that intersect $Q_i^{k*}$. In fact, observe that if
$Q_i^{k*}\bigcap Q_j^{k+1*}\not=\o$,  by Lemma 5.1, we have
 $$Q_j^{k+1*}\subset
2^6nQ_j^{k*}\subset \oz^k.$$ So,  we set $\wt
Q_i^k:=2^6nQ_j^{k*}.$ Note that $l_j^{k+1}<1$ and $l_j^k<1$,
then
$$ \wt Q_i^k\subset \oz^k,\ {\rm and}\ \wz(\wt Q_i^k)\le C_{10}\wz(Q_i^{k*}).$$
Moreover, $h_i^k$ satisfies the moment conditions. This is clear
for $(f-P_i^k)\eta_i^k$ and
$(f-P_j^{k+1})\eta_j^{k+1}\eta_i^k+P_{ij}^{k+1}\eta_i^{k+1}$.

Let $\lz_i^k=C_{10}2^k[\wz(\wt Q_i^k)]^{1/p}$ and
$a_i^k=(\lz_i^k)^{-1}h_i^k$. Moreover, by (5.3) and above Cases I,
II and III, we  know that $a_i^k$ is a $(p,\fz,s)^\gz_\wz$-atom.
By $\wz\in A_q^{loc}$ and Proposition 2.1(i), we have
$$\begin{array}{cl}
\dsum_{k\in\zz}\dsum_{i\in\nn}|\lz_i^k|^p&\le
C\dsum_{k\in\zz}\dsum_{i\in\nn}2^{kp}\wz(\wt Q_i^k)\le
C\dsum_{k\in\zz}\dsum_{i\in\nn}2^{kp}\wz( Q_i^{k^*})\\
&\le C\dsum_{k\in\zz}2^{kp}\wz(\oz_k)\le C\|{\cal
M}_N(f)\|_{L^p_\wz(\rz)}^p\le C\|f\|_{h^p_{\wz,N}(\rz)}^p.
\end{array}$$
We now consider the case $k_0>-\fz$, which together with $f\in
h^p_{\wz,N}(\rz)$ implies $\wz(\rz)<\fz$. Adapting the previous
arguments, we have
$$f=\dsum_{k=k_0}^\fz(g^{k+1}-g^k)+g^{k_0}:=\wt f+g^{k_0}.$$
For the function $\wt f$, we have the same $(p,\fz,s)_\wz$
atomic decomposition as above and
$$\dsum_{k\ge k_0}\dsum_{i\in\nn}|\lz_i^k|^p\le
C\| f\|_{h^p_{\wz,N}(\rz)}^p.$$  For the function $g^{k_0}$, it is easy to see that there
exists a positive constant $C_{11}$ such that
$$\|g^{k_0}\|_{L^\fz_\wz(\rz)}\le C_{11} 2^{k_0}\le
2C_{11}\dinf_{x\in\rz}{\cal M}_Nf(x).$$ Let
$$a_0(x)=g^{k_0}(x)2^{-k_0}C_{11}^{-1}[\wz(\rz)]^{-1/p}, \quad \lz_0=C_{11}2^{k_0}[\wz(\rz)]^{1/p}.$$
Hence,
 $$|\lz_0|^p\le (2C_{11})^p\| f\|_{h^p_{\wz,N}(\rz)}^p, \ {\rm and}\ \|a_0\|_{L^\fz_\wz(\rz)}\le [\wz(\rz)]^{-1/p}.$$
Then,
$$\dsum_{k\ge k_0}\dsum_{i\in\nn}|\lz_i^k|^p
+|\lz_0|^p\le C\| f\|_{h^p_{\wz,N}(\rz)}^p.$$ The  proof of Lemma
5.4 is complete.

{\bf Remark 5.1:}\quad In fact, from the proof of Lemma 5.4, we
can take all  $(p,\fz,s)_\wz$ atoms with sidelengths $\le 2$ in Lemma 5.4.

The following is one of the main results in this paper.
\begin{thm}\label{t5.1.}\hspace{-0.1cm}{\rm\bf 5.1.}\quad
Let $\wz\in A_\fz^{loc}$ and $q_\wz$ be as in (2.3). If $q\in
(q_\wz,\fz], p\in (0,1]$, $N\ge N_{p,\wz}$, and
$s\ge[n(q_\wz/p-1]$, then $h^{p,q,s}_{\wz}(\rz)=h^p_{\wz, N}(\rz)=
h^p_{\wz, N_{p,\wz}}(\rz)$ with equivalent norms.
\end{thm}
Proof. It is easy to see that
$$h^{p,\fz,\bar s}_{\wz}(\rz)\subset h^{p,q,s}_{\wz}(\rz)\subset
h^p_{\wz, N_{p,\wz}}(\rz)\subset h^p_{\wz, N}(\rz) \subset
h^p_{\wz, \bar N}(\rz),$$ where $\bar s$ is an integer no less
than $s$ and $\bar N$ is an integer larger than $N$, and the
inclusions are continuous. Thus, to prove Theorem 5.1, it suffices
to prove that for any $N>s\ge[q_\wz/n]$, $h^p_{\wz, N}(\rz)\subset
h^{p,\fz,s}_{\wz}(\rz)$, and for all $f\in h^p_{\wz, N}(\rz)$,
$\|f\|_{h^{p,\fz,s}_{\wz}(\rz)}\le C\|f\|_{h^p_{\wz, N}(\rz)}$.

To this end, let $f\in h^p_{\wz,N}(\rz)$. By Corollary 4.1, there
exists a sequence of functions, $\{f_m\}_{m\in\nn}\subset
(h^p_{\wz,N}(\rz)\bigcap L_\wz^1(\rz))$, such that
$\|f_m\|_{h^p_{\wz, N}(\rz)}\le 2^{-m}\le \|f\|_{h^p_{\wz,
N}(\rz)}$ and $f=\dsum_{m\in\nn}f_m$ in $h^p_{\wz, N}(\rz)$. By
Lemma 5.4, for each $m\in\nn$, $f_m$ has an atomic decomposition
$f_m=\sum_{i\in\nn_0}\lz_i^ma_i^m$ in ${\cal D}'(\rz)$, where
$\sum_{i\in\nn_0}|\lz_i^m|^p\le C\|f_m\|_{h^p_{\wz,N}(\rz)}^p$ and
$\{a_i^m\}_{i\in\nn_0}$ in ${\cal D}'(\rz)$, where
$\sum_{i\in\nn_0}|\lz_i^m|^p\le C\|f_m\|_{h^p_{\wz,N}(\rz)}^p$ and
$\{a_i^m\}_{i\in\nn_0}$ are $(p,\fz,s)_\wz$-atoms. Since
$$\dsum_{m\in\nn_0}\dsum_{i\in\nn_0}|\lz_i^m|^p\le C\dsum_{m\in\nn_0}
\|f_m\|_{h^p_{\wz,N}(\rz)}^p\le C\|f\|_{h^p_{\wz,N}(\rz)}^p,$$
then $f=\sum_{m\in\nn_0}\sum_{i\in\nn_0}\lz_i^ma_i^m\in
h^{p,\fz,s}_{\wz}(\rz)$ and $\|f\|_{h^{p,\fz,s}_{\wz}(\rz)}\le
C\|f\|_{h^p_{\wz,N}(\rz)}$. Thus, Theorem 5.1 is proved.

For simplicity, from now on, we denote by $h^p_{\wz}(\rz)$ the
weighted local Hardy space $h^p_{\wz, N}(\rz)$ associated with
$\wz$, where $N\ge N_{p,\wz}$. Moreover, it is easy to see that
$h^1_\wz\subset L_\wz^1(\rz)$ via weighted atomic decomposition.
However, the elements in $h^p_\wz(\rz)$ with $p(0,1)$ are not
necessary functions thus $h^p_\wz(\rz)\not= L_\wz^p(\rz)$. But,
for any $q\in (q_\wz,\fz)$, by Lemma 5.4 and pointwise
convergence of weighted atomic decompositions, we have
$(h^p_\wz(\rz)\bigcap L_\wz^1(\rz))\subset L_\wz^p(\rz)$,  and
for all $f\in (h^p_\wz(\rz)\bigcap L_\wz^1(\rz))$,
$\|f\|_{L^p_\wz(\rz)}\le\|f\|_{h^p_\wz(\rz)}$.

\begin{center} {\bf6. Finite atomic decompositions}\end{center}
In this section, we prove that for any given finite linear
combination of weighted atoms when $q<\fz$, its norm in
$h^p_\wz(\rz)$ can be achieved via all its finite weighted atomic
decompositions. This extends the main results in \cite{msv} to the
setting of weighted local Hardy spaces.

Let $\wz\in A_\fz^{loc}$ and $(p,q,s)_\wz$ be an admissible
triplet. Denote by $h^{p,q,s}_{\wz,fin}(\rz)$ the vector space of
all finite linear combination of $(p,q,s)_\wz$-atoms and single
atom, and the norm of $f$ in $h^{p,q,s}_{\wz,fin}(\rz)$ is defined
by
$$\begin{array}{cl}
\|f\|_{h^{p,q,s}_{\wz,fin}(\rz)}&=\dinf\l\{\l[\dsum_{j=0}^k|\lz_j|^p\r]^{1/p}:
f=\dsum_{j=0}^k\lz_ja_j, k\in\nn_0, \{a_i\}_{i=1}^k\ {\rm are}\ (p,q,s)_\wz-\r. \\
&\qquad\l. {\rm atoms\ with\ sidelenghths\ }\le 2, {\rm and}\ a_0\ {\rm is\ a}\ (p,q)_\wz\ {\rm single\
atom}\r\}.\end{array}$$ Obviously, for any admissible triplet
$(p,q,s)_\wz$ atom and $(p,q)_\wz$ single atom, the set
$h_{\wz,fin}^{p,q,s}(\rz)$ is dense in $h^{p,q,s}_{\wz}(\rz)$ with
respect to the quasi-norm $\|\cdot\|_{h_{\wz, fin}^{p,q,s}(\rz)}$.
\begin{thm}\label{t6.1.}\hspace{-0.1cm}{\rm\bf 6.1.}\quad
Let $\wz\in A_\fz^{loc}$, $q_\wz$ be as in (2.3), and
$(p,q,s)_\wz$ be an admissible triplet with sidelenghth $\le 2$. If $q\in (q_\wz,\fz)$,
then $\|\cdot\|_{h_{\wz,fin}^{p,q,s}(\rz)}$ and
$\|\cdot\|_{h_{\wz}^p(\rz)}$ are equivalent quasi-norms on
$h_{\wz,fin}^{p,q,s}(\rz)$.
\end{thm}
Proof. Clearly, $\|f\|_{h^p_\wz(\rz)}\le
\|f\|_{h^{p,q,s}_{\wz,fin}(\rz)}$ for $f\in
h^{p,q,s}_{\wz,fin}(\rz)$ and for $q\in (q_\wz,\fz)$. Thus,we have
to show that for every $q$ in $(q_\wz,\fz)$ there exists a
constant $C$ such that for all $f\in h^{p,q,s}_{\wz,fin}(\rz)$
$$\|f\|_{h^{p,q,s}_{\wz,fin}(\rz)}\le C\|f\|_{h^p_\wz(\rz)}.\eqno(6.1)$$
Suppose that $q\in (q_\wz,\fz)$ and that $f$ is in
$h^{p,q,s}_{\wz,\gz,fin}(\rz)$ with $\|f\|_{h^p_{\wz}(\rz)}=1$. In
this section, we take $k_0\in\zz$ such that
$2^{k_0-1}\le\dinf_{x\in\rz}{\cal M}_Nf(x)<2^{k_0}$, if
$\dinf_{x\in\rz}{\cal M}_Nf(x)=0$, write $k_0=-\fz$.  For each
integer $k\ge k_0$, set
$$\oz_k\equiv\{x\in\rz: {\cal M}_Nf(x)>2^k\},$$
where and in what follows $N=N_{p,\wz}$. We use the same notation
as in Lemma 5.4. We first consider the case $k_0=-\fz$. Since $f\in (h_\wz^p(\rz)\bigcap L_\wz^q(\rz))$,
by Lemma 5.4, there exists numbers $\{\lz_i^k\}_{k\in\zz,
i\in\nn}\subset \cc$ and $(p,\fz,s)_\wz$-atoms $\{a_i^k\}_{k,i\in\nn}$, $\lz_0\subset \cc$  such that
$$f= \sum_{k}\sum_{i\in\nn}\lz_i^ka_i^k$$ holds almost everywhere and
in ${\cal D}'(\rz)$, and (i) and (ii) in Lemma 5.4 hold.

Obviously, $ f$ has compact support.
 Suppose that $\supp  f\subset
Q(x_0,r_0)$. We write $\bar Q=Q(x_0, 2^{3(10+n)}r_0+2n)$. For $\vz$
in ${\cal D}_N$ and $x\in \rz\setminus \bar Q$, for $0<t<1$, we
have
$$\vz_t* f(x)=0.$$
Hence, $\supp \sum_{k}\sum_{i\in\nn}\lz_i^ka_i^k\subset
\bar Q$.

We claim that the series $ \sum_{k}\sum_{i\in\nn}\lz_i^ka_i^k$ converges to $f$ in
$L^q_\wz(\rz)$. For any $x\in\rz$, since
$\rz=\bigcup_{k\in\zz}(\oz_k\setminus\oz_{k+1})$, there exists
$j\in\zz$ such that $x\in (\oz_j\setminus\oz_{j+1})$. Since $\supp
a_i^k\subset Q_i^k\subset \oz_k\subset\oz_{j+1}$ for $k>j$, then
applying Lemmas 5.1, 5.2 and  Lemma 5.4, we have
$$|\sum_{k}\sum_{i\in\nn}\lz_i^ka_i^k|\le C\dsum_{ k\le j}2^k
\le C2^j\le C{\cal M}_Nf(x).$$ Since $f\in L^q_\wz(\rz)$, we have
${\cal M}_Nf\in L^q_\wz(\rz)$. The Lebesgue dominated convergence
theorem now implies that $\sum_{k}\sum_{i\in\nn}\lz_i^ka_i^k$ converges to $f$ in
$L^q_\wz(\rz)$, and the claim is proved.

For each positive integer $K$ we denote by $F_K=\{(i,k):
k, |i|+|k|\le K\}$ and $f_K=\sum_{(i,k)\in
F_K}\lz_i^ka_i^k$. Observing that for any $\ez\in (0,1)$,
if $K$ is large enough, by $ f\in L^q_\wz$, we have $ (f-f_K)/\ez$
is a $(p,q,s)_\wz$-atom. Since $(f-f_K)/\ez\in\bar Q=Q(x_0, 2^{3(10+n)}r_0+2n)$, so we
can divide $\bar Q$  into  $N_0$ (depending only on $r_0$ and $n$) disjoint cubes $\{Q_i\}_{i=1}^{N_0}$
with sidelengths $1\le l_i\le 2$. Then, $ (f-f_K)\chi_{Q_i}/\ez$
is a $(p,q,s)_\wz$-atom for $i=1,\cdots, N_0$.
Thus, $f=f_K+\sum_{i=1}^{N_0}(f-f_K)\chi_{Q_i}$ is a  linear
weighted atom combination of $f$. Taking $\ez=N_0^{-1/p}$ and by Lemma 5.4, we have
$$\|f\|_{h^{p,q,s}_{\wz,fin}(\rz)}\le \dsum_{(i,k)\in F_K}|\lz_i^k|^p+N_0\ez^p\le C.$$

We now  consider the case $k_0>-\fz$. Since $f\in (h_\wz^p(\rz)\bigcap L_\wz^q(\rz))$,
by Lemma 5.4, there exists numbers $\{\lz_i^k\}_{k\in\zz,
i\in\nn}\subset \cc$ and $(p,\fz,s)_\wz$-atoms $\{a_i^k\}_{k\ge
k_0,i\in\nn}$, $\lz_0\subset \cc$ and the $(p,\fz)_\wz$ singe atom
$a_0$ such that
$$f= \sum_{k\ge k_0}\sum_{i\in\nn}\lz_i^ka_i^k+\lz_0a_0$$ holds almost everywhere and
in ${\cal D}'(\rz)$, and (i) and (ii) in Lemma 5.4 hold.
As the case $k_0=-\fz$, we can prove that the series $ \sum_{k\ge
k_0}\sum_{i\in\nn}\lz_i^ka_i^k+\lz_0a_0$ converges to $f$ in
$L^q_\wz(\rz)$.

Finally, for each positive integer $K$ we denote by $F_K=\{(i,k):
k\ge k_0, |i|+|k|\le K\}$ and $f_K=\sum_{(i,k)\in
F_K}\lz_i^ka_i^k+\lz_0a_0$. If $K$ is large enough, then $\|f-f_K\|_{L^q(\wz)}\le
[\wz(\rz)]^{1/q-1/p}$. So, $(f-f_K)$ is a $(p,q)_\wz$ single atom.
 By Lemma 5.4, we have
$$\|f\|_{h^{p,q,s}_{\wz,fin}(\rz)}\le \dsum_{(i,k)\in F_K}|\lz_i^k|^p+\lz_0^p\le C.$$
Thus, (6.1) holds. The proof is finished.

 As an application of finite atomic
decompositions, we establish boundedness in $h^p_\wz(\rz)$ of
quasi- Banach-valued sublinear operators.

As in \cite{blyz}, we recall that a quasi-Banach space ${\cal B}$
is a vector space endowed with a quasi-norm $\|\cdot\|_{{\cal B}}$
which is nonnegative, non-degenerate (i.e., $\|f\|_{\cal B}=0$ if
and only if $f=0$), homogeneous, and obeys the quasi-triangle
inequality, i.e., there exists a positive constant $K$ no less
than $1$ such that for all $f,g\in{\cal B}$, $\|f+g\|_{{\cal
B}}\le K(\|f\|_{\cal B}+\|g\|_{\cal B})$.

Let $\bz\in (0,1]$. A quasi-Banach space ${\cal B}_\bz$ with the
quasi-norm $\|\cdot\|_{{\cal B}_\bz}$ is said to be a
$\bz$-quasi-Banach space if $\|f+g\|^\bz_{{\cal B}_\bz}\le
\|f\|^\bz_{{\cal B}_\bz}+\|g\|^\bz_{{\cal B}_\bz}$ for all $f,g\in
{\cal B}_\bz$.

Notice that any Banach space is a $1$-quasi-Banach space, and the
quasi-Banach space $l^\bz,\ L^\bz_\wz(\rz)$ and $h^\bz_\wz(\rz)$
with $\bz\in (0,1)$ are typical $\bz$-quasi-Banach spaces.

For any given $\bz$-quasi-Banach space ${\cal B}_\bz$ with $\bz\in
(0,1]$ and  a linear space ${\cal Y}$, an operator $T$ from
 ${\cal Y}$ to ${\cal B}_\bz$ is said to be ${\cal B}_\bz$-sublinear
 if for any $f, g\in {\cal B}_\bz$  and $\lz,\ \nu\in\cc$, we have
 $$\|T(\lz f+\nu g)\|_{{\cal B}_\bz}\le\l(|\lz|^\bz\|T(f)\|_{{\cal
 B}_\bz}^\bz
+|\nu|^\bz\|T(g)\|_{{\cal B}_\bz}^\bz\r)^{1/\bz}$$ and
$\|T(f)-T(g)\|_{{\cal B}_\bz}\le \|T(f-g)\|_{{\cal B}_\bz}$.

We remark that if $T$ is linear, then $T$ is ${\cal
B}_\bz$-sublinear. Moreover, if ${\cal B}_\bz=L_\wz^q(\rz)$, and
$T$ is nonnegative and sublinear in the classical sense, then $T$
is also ${\cal B}_\bz$-sublinear.
\begin{thm}\label{t6.2.}\hspace{-0.1cm}{\rm\bf 6.2.}\quad
Let $\wz\in A_\fz^{loc}, 0<p\le\bz\le 1$, and ${\cal B}_\bz$ be a
$\bz$-quasi-Banach space. Suppose $q\in (q_\wz,\fz)$ and $T:
h^{p,q,s}_{\wz,fin}(\rz)\to {\cal B}_\bz$ is a ${\cal
B}_\bz$-sublinear operator  such that
$$\begin{array}{cl}S\equiv\{\|T(a)\|_{{\cal B}_\bz}:&\ a \ {\rm is\ any\ }\
(p,q,s)_\wz-{\rm atom\ with\ sidelength\ }\le 2  \\
&\quad {\rm or}\ (p,q)_\wz\ {\rm  single\
atom}\}<\fz.\end{array}$$ Then there exists a unique bounded ${\cal
B}_\bz$-sublinear operator $\wt T$ from $h^p_\wz(\rz)$ to ${\cal
B}_\bz$ which extends $T$.
\end{thm}
Proof. For any $f\in h^{p,q,s}_{\wz,fin}(\rz)$, by Theorem 6.1,
there exist numbers $\{\lz_j\}_{j=0}^l\subset\cc$ and
$(p,q,s)_\wz$-atoms $\{a_j\}_{j=1}^l$ and  the $(p,q)_\wz$ single
atom $a_0$ such that $f=\sum_{j=0}^l\lz_ja_j$ pointwise and
$\sum_{j=0}^l|\lz_j|^p\le C\|f\|_{h^p_\wz(\rz)}^p$. Then by the
assumption, we have
$$\|T(f)\|_{{\cal B}_\bz}\le C\l[\dsum_{j=0}^l|\lz_j|^p\r]^{1/p}\le
C\|f\|_{h^p_\wz(\rz)}.$$ Since $h^{p,q,s}_{\wz,fin}(\rz)$ is dense
in $h^p_\wz(\rz)$, a density argument gives the desired results.

\begin{center} {\bf 7. Applications}\end{center}
In this section, we study weighted $L^p$ inequalities for strongly
singular integrals  and pseudodifferential operators and their
commutators.

Given a  real number $\tz>0$ and a smooth radial cut-off function
$v(x)$ supported in the ball $\{x\in\rz: \ |x|\le 2\}$, we
consider the strongly singular kernel
$$k(x)=\dfrac{e^{i|x|^{-\tz}}}{|x|^n}v(x).$$
Let us denote by $Tf$ the corresponding strongly singular integral
operator:
$$Tf(x)=p.v\dint_\rz k(x-y)f(y)dy.$$
This operator has been studied by several authors, see \cite{h},
\cite{w}, \cite {fs}, \cite{c} and \cite{ghst}. In particular, S.
Chanillo \cite{c} established  the weighted $L^p_\wz(\rz)$
boundedness for strongly singular integrals provided that $\wz\in
A_p(\rz)$ (Muckenhoupt weights) for $1<p<\fz$. J. Garc\'ia-Cuerva
et al \cite{ghst} obtained weighted $L^p$ estimates with pairs of
weights for commutators generated by the strongly singular
integrals and the classical $BMO(\rz)$ functions. We have the
following results for the strongly singular integrals.
\begin{thm}\label{t7.1.}\hspace{-0.1cm}{\rm\bf 7.1.}\quad
Let $T$ be strongly singular integral operators, then
\begin{enumerate}
\item[(i)]$\|Tf\|_{L^p_\wz(\rz)}\le C_{p,\wz}\|f\|_{L^p_\wz(\rz)}$
for $1<p<\fz$ and $\wz\in A_p^{loc}$.
\item[(ii)]$\|Tf\|_{L^{1,\fz}_\wz(\rz)}\le
C_\wz\|f\|_{L^1_\wz(\rz)}$ for $\wz\in A_1^{loc}$.
\item[(iii)]$\|Tf\|_{L^1_\wz(\rz)}\le C_\wz\|f\|_{h^1_\wz(\rz)}$
for $\wz\in A_1^{loc}$.
\end{enumerate}
\end{thm}
Proof.  We first note that for $\wz\in A_p$ the inequality (i) is
known to be true, see \cite{c}. For $\wz\in A_p^{loc}$, by Lemma
2.1 (i) for any unit cube $Q$ there is a $\bar\wz\in A_p$ so that
$\bar\wz=\wz$ on $6Q$. Then
$$\begin{array}{cl}
\|Tf\|_{L^p_\wz(Q)}&= \|T(\chi_{6Q}f)\|_{L^p_\wz(Q)}\\
&\le \|T(\chi_{6Q}f)\|_{L^p_{\bar\wz}(Q)}\\
&\le C\|(\chi_{6Q}f)\|_{L^p_{\bar\wz}(\rz)}\\
&\le C\|f\|_{L^p_{\bar\wz}(6Q)}.
\end{array}$$
Summing over all dyadic unit $I$ gives (i).

For (ii), similar to (i), note that for $\wz\in A_1$ the
inequality (ii) is known to be true, see \cite{c}. Since $\wz\in
A_p^{loc}$, by Lemma 2.1 (i) for any unit cube $I$ there is a
$\bar\wz\in A_1$ so that $\bar\wz=\wz$ on $6Q$. Then for any
$\lz>0$
$$\begin{array}{cl}
\wz(\{x\in Q: |Tf(x)|>\lz\})&\le \wz(\{x\in Q: |T(\chi_{6Q}f)(x)|>\lz\})\\
&= \bar\wz(\{x\in Q: |T(\chi_{6Q}f)(x)|>\lz\})\\
&\le C\lz^{-1}\|(\chi_{6Q}f)\|_{L^1_{\bar\wz}(\rz)}\\
&= C\lz^{-1}\|f\|_{L^1_{\wz}(6Q)}.
\end{array}$$
Summing over all dyadic unit $Q$ gives (ii).

Finally, to consider (iii). Let $a(x)$ be an atom in
$h^1_\wz(\rz)$, supported in a cube $Q$ centered at $x_0$ and
sidelength $\dz\le 2$ by Remark 5.1, or $a(x)$ is a single atom.
To prove the (iii), by Theorem 6.2, it is enough to show that
$$\|Ta\|_{L^1_\wz(\rz)}\le C,\eqno(7.1)$$
where $C$ is independent of $a$.

 It is easy to see that (7.1) holds while $a(x)$ is a
single atom. It remains to consider this kind of atom supported in
a cube $Q$ centered at $x_0$ and sidelength $\dz$.
 In deed, let $\dz_0$ be a number satisfying
$4\dz_0=\dz_0^{1/(1+\tz)}$. Obviously, $\dz_0<1$.

Case 1. $2\ge\dz\ge  \dz_0$. This is the trivial case. Let
$Q^*=(10n/\dz_0)Q$. Now
$$\dint_\rz |Ta|\wz(x)dx=\dint_{Q^*} |Ta|\wz(x)dx
+\dint_{\rz\setminus Q^*} |Ta|\wz(x)dx=\dint_{Q^*} |Ta|\wz(x)dx.$$
Obviously,
$$\begin{array}{cl}
\dint_{Q^*} |Ta|\wz(x)dx&\le
C\l(\dint_{\rz}|Ta|^p\wz(x)dx\r)^{1/p}\l(\dint_{Q^*}\wz(x)dx\r)^{1/p'}\\
&\le
C\l(\dint_{\rz}|a|^p\wz(x)dx\r)^{1/p}\l(\dint_{Q^*}\wz(x)dx\r)^{1/p'}\\
&\le C\wz(Q)^{-1/p'}\l(\dint_{Q^*}\wz(x)dx\r)^{1/p'}\le C.
\end{array}\eqno(7.2)$$
Case 2. $\dz< \dz_0$. We let $Q^*=4Q$ and $\bar Q=Q(x_0,\dz^{1/(1+\tz)})$.
Then
$$\begin{array}{cl}
\dint_\rz|T a|\wz(x)dx&\le \dint_{Q^*}|T a|\wz(x)dx+ \dint_{\bar
Q\setminus Q^*}|T a|\wz(x)dx+\dint_{\rz\setminus \bar Q}|T
a|\wz(x)dx\\
&:=I+II+III.\end{array}$$ For $I$, similar to (7.2), we have
$$I\le
C\l(\dint_{\rz}|Ta|^p\wz(x)dx\r)^{1/p}\l(\dint_{Q^*}\wz(x)dx\r)^{1/p'}\le
C.$$ We now estimate the term $III$. Clearly, by the mean value
theorem,
$$\begin{array}{cl}
|Ta(x)|&\le \dfrac {C\dz}{|x-x_0|^{\tz+n+1}}\chi_{\{|x-x_0|<4n\}}(x)\dint_Q|a(y)|dy\\
&\le \dfrac {C\dz}{|x-x_0|^{\tz+n+1}}\chi_{\{|x-x_0|<4n\}}(x)\\
&\qquad\times\l(\dint_Q|a(x)|^p\wz(x)dx\r)^{1/p}\l(\dint_{Q}[\wz(x)]^{-p'/p}dx\r)^{1/p'}\\
&\le \dfrac
{C\dz}{|x-x_0|^{\tz+n+1}}\chi_{\{|x-x_0|<4n\}}(x)\dfrac
{|Q|}{\wz(Q)}.
\end{array}$$
Hence, by the properties of $A_1^{loc}$ (see Lemma 2.1), we have
$$\begin{array}{cl}
III&\le \dfrac
{C\dz|Q|}{\wz(Q)}\dint_{\dz^{1/(1+\tz)}\le|x-x_0|\le
4n}\dfrac {\wz(x)}{|x-x_0|^{\tz+n+1}}dx\\
&\le \dfrac {C\dz|Q|}{\wz(Q)}\dsum_{k=k_0}^{k_1}\dfrac
1{(2^k\dz)^{1+\tz}}\l(\dfrac 1{(2^k\dz)^n}\dint_{|x-x_0|\le
2^k\dz}\wz(x)dx\r)\\
&\le C,
\end{array}$$
where $k_0$  and $k_1$ are positive integers such that
$2^{k_0}\dz\le \dz^{1/(1+\tz)}\le 2^{k_0+1}\dz$ and $2^{k_1-1}\le
4n\le 2^{k_1}$. We now estimate the term $II$. For $x\in \bar
Q\setminus Q^*$
$$\begin{array}{cl}
Ta(x) &=\dint_\rz\dfrac
{e^{i|x-y|^{-\tz}}v(x-y)}{|x-y|^{n(2+\tz)/r'}}\\
&\qquad\times\l(\dfrac1{|x-y|^{n(1-(2+\tz)/r')}}
-\dfrac1{|x_0-x|^{n(1-(2+\tz)/r')}}\r)a(y)dy\\
&\qquad+\dint_\rz\dfrac
{e^{i|x-y|^{-\tz}}v(x-y)}{|x-y|^{n(2+\tz)/r'}}\dfrac{a(y)}{|x_0-x|^{n(1-(2+\tz)/r')}}dy\\
&=A(x)+B(x),
\end{array}$$
where $r'$ is taken so close to 1 to guarantee that $2+\tz<r$.
Applying the mean value theorem to the term in brackets in the
integrand of $A$, and noting that for $y\in Q$, and $x\in \bar
Q\setminus Q^*$, $|x-y|\ge c|x-x_0|$, we have

$$|A(x)|\le
\dfrac
{C\dz}{|x-x_0|^{n+1}}\chi_{\{|x-x_0|<4n\}}(x)\dint_Q|a(y)|dy\le
\dfrac {C|Q|}{|x-x_0|^{n+1}}\chi_{\{|x-x_0|<4n\}}(x)\dfrac
{|Q|}{\wz(Q)}.$$ Therefore,
$$\begin{array}{cl}
II&\le \dfrac {C\dz|Q|}{\wz(Q)} \dint_{\dz\le|x-x_0|\le
4n}\dfrac {\wz(x)}{|x-x_0|^{n+1}}dx\\
&\qquad +C\dint_{\dz\le|x-x_0|<\dz^{1/(1+\tz)}}|K_{\tz,r}*a|\dfrac{\wz(x)}{|x_0-x|^{n(1-(2+\tz)/r')}}dx\\
&\le C+C\l(\dint_\rz
|K_{\tz,r}*a|^rdx\r)^{1/r}\l(\dint_{\dz<|x-x_0|<\dz^{1/(1+\tz)}}
\dfrac{\wz(y)^{r'}}{|x_0-x|^{n(1-(2+\tz)/r')}}dx\r)^{1/r'}\\
&\le
C+C\|a\|_{L^{r'}(\rz)}\l(\dsum_{k=0}^{k_0}(2^k\dz)^{(r'-1)(\tz+1)}\dfrac
1{(2^k\dz)^n}\dint_{|x-x_0|\le 2^k\dz}\wz(x)^{r'}dx\r)^{1/r'}\\
&\le C,
\end{array}$$
where $2^{k_0-1}\dz<\dz^{1/(1+\tz)}\le 2^{k_0}\dz$,
$K_{\tz,r}(x):=\frac{e^{i|x|^{-\tz}}}{|x|^{(\tz+2)/r}}$, and we
used the following fact (see \cite{c})
$$\|K_{\tz,r}*f\|_{L^r(\rz)}\le C_r\|f\|_{L^{r'}(\rz)},\quad
r>2+\tz.$$ Thus,  Theorem 7.1 is proved.

We now introduce $BMO_{loc}$ of locally integrable functions with
bounded mean oscillation which has a  intimate relationship
between the $A_p^\loc$ weights.  Namely,
$$\|b\|_{BMO^{loc}}:=\dsup_{|Q|\le 1}\dfrac 1{|Q|}\dint_Q|b-b_Q|\,dx<\fz,$$
where $b_Q=\frac 1{|Q|}\int_Q f(x)\,dx$.

It is easy to see that, we have the following result.
\begin{lem}\label{l7.1.}\hspace{-0.1cm}{\rm\bf 7.1.}\quad
  Fix $p>1$ and let $b\in BMO^{loc}$. Then there
exists $\ez>0$, depending upon the $BMO^{loc}$ constant of $b$,
such that $e^{xb}\in A_p^{loc}$ for $|x|<\ez$.
\end{lem}
\begin{lem}\label{l7.2.}\hspace{-0.1cm}{\rm\bf 7.2.}\quad
 Let $b\in BMO^{loc}$, then there exist  positive constants $c_1$ and $c_2$ such that
for every cube $Q$ with $|Q|\le 1$ and every $\lz>0$, we have
$$|\{x\in Q: |b(x)-b_Q|>\lz\}|\le c_1|Q|\exp\l\{-\frac
{c_2\lz}{\|b\|_{BMO^{loc}(\rz)}}\r\}.$$
\end{lem}
As a consequence of Lemma 7.2 and Lemma 2.1, we have the following
result.
\begin{cor}\label{c7.1.}\hspace{-0.1cm}{\rm\bf 7.1.}\quad
 Let $b\in BMO^{loc}$ and $\wz\in A_\fz^{loc}$, then there exist  positive constants $C_3$ and $C_4$ such that
for every cube $Q$ with $|Q|\le 1$ and every $\lz>0$, we have
$$\wz(\{x\in Q: |b(x)-b_Q|>\lz\})\le c_3\wz(Q)\exp\l\{-\frac
{c_4\lz}{\|b\|_{BMO^{loc}(\rz)}}\r\}.$$
\end{cor}
As an application of Corollary 7.2, we have
\begin{prop}\label{p7.1.}\hspace{-0.1cm}{\rm\bf 7.1.}\quad
 Let $b\in BMO^{loc}$, $1\le p<\fz$, and $\wz\in A_\fz^{loc}$, then there exists  a positive constant $C$  such that
for every cube $Q$ with $|Q|\le 1$
$$\dfrac 1{\wz(Q)}\dint_Q|b(x)-b_Q|^p\wz(x)dx\le C\|b\|_{BMO^{loc}}^p.$$
\end{prop}

 We now consider in this paper commutator
of Coifman-Rochberg-Weiss $[b,T]$ defined by the formula
$$[b,T]f(x)=b(x)Tf(x)-T(bf)(x)=\dint_\rz(b(x)-b(y))k(x-y)f(y)dy.$$
 As in the case of strongly singular integrals, we have

\begin{thm}\label{t7.2.}\hspace{-0.1cm}{\rm\bf 7.2.}\quad
Let  $b\in BMO_{loc}(\rz)$ and $T$ be the strongly singular integral
operators , then
\begin{enumerate}
\item[(i)]$\|[b,T]f\|_{L^p_\wz(\rz)}\le C_{p,\wz}\|b\|_{BMO_{loc}(\rz)}\|f\|_{L^p_\wz(\rz)}$ for
$1<p<\fz$ and $\wz\in A_p^{loc}$.
\item[(ii)]$\|[b,T]f\|_{L^{1,\fz}_\wz(\rz)}\le C_\wz\|b\|_{BMO_{loc}(\rz)}\|f\|_{h^1_\wz(\rz)}$ for $\wz\in
A_1^{loc}$.

\end{enumerate}
\end{thm}
 Proof:\quad By  of Lemma 7.1, there is $\eta>0$ such that
$\wz^{(1+\eta)}\in A_p^{loc}$. Then, we choose $\dz>0$ such that
$exp(s\dz b(1+\eta)/\eta)\in A_p^{loc}$ if $0\le
s(1+\eta)/\eta<\dz$ with uniform constant. For $z\in {\cal C}$ we
define the operator
$$T_zf=e^{zb}T(e^{-zb}f).$$
We claim that
$$\|T_z f\|_{L^p_\wz(\rz)}\le C\|f\|_{L^p_\wz(\rz)}$$
uniformly on $|z|\le s<\dz\eta/(1+\eta)$.

The function $z\to T_z f$ is analytic, and by the Cauchy theorem,
if $s<\dz\eta/(1+\eta)$,
$$\dfrac d{dz} T_z f|_{z=0}=\dfrac 1{2\pi i}\dint_{|z|=s}\dfrac {T_z
f}{s^2}\,dz.$$ Observing that
$$\dfrac d{dz} T_z f|_{z=0}=[b, T]f$$
and applying the Minkowski inequality to the previous equality, we
get
$$\|[b, T]f\|_{L^p_\wz(\rz)}\le  \dfrac 1{2\pi}\dint_{|z|=s}\dfrac {\|T_z
f\|_{L^p(\wz)}}{s^2}\,|dz|\le \dfrac Cs \|f\|_{L^p_\wz(\rz)}.$$ It
remains to prove the claim, which is equivalent to
$$\begin{array}{cl}
&\l(\dint_\rz |T f(x)|^p exp({\cal
R}(z)qb(x))\wz(x)\,dx\r)^{1/p}\\
&\qquad\qquad\qquad\le C \l(\dint_\rz |f(x)|^pexp({\cal
R}(z)pb(x))\wz(x)\,dx\r)^{1/p}.\end{array}\eqno(7.3)$$ We write
$\wz_0:=exp({\cal R}(z)b(1+\eta)/\eta)$ and $\wz_1:=
\wz^{1+\eta}$. Since $\wz_0$ and $\wz_1\in A_p^{\loc}$, we have
$$\l(\dint_\rz |T f(x)|^p\wz_0(x)\,dx\r)^{1/p}\le C
\l(\dint_\rz |f(x)|^p\wz_0(x)\,dx\r)^{1/p}$$ and
$$\l(\dint_\rz |T f(x)|^p\wz_1(x)\,dx\r)^{1/p}\le C
\l(\dint_\rz |f(x)|^p\wz_1(x)\,dx\r)^{1/p}.$$ Now, by Stein-Weiss
interpolation theorem, we have
$$\l(\dint_\rz |T f(x)|^p\wz_0^{(1-\bz)}\wz_1^{\bz}\,dx\r)^{1/p}\le C
\l(\dint_\rz |f(x)|^p\wz_0^{(1-\bz)}\wz_1^{\bz}\,dx\r)^{1/p}$$ and
taking $\bz=(1+\eta)^{-1}$, then we obtain (7.3). Thus, (i) of
Theorem 7.2 is proved.

For (ii), Let the function $a_j(x)$ is a $h^1_\wz(\rz)$ atom  and
$\supp\ a_j\subset Q(x_j,r_j)$, and $a_0$ is a single atom if
$\wz(\rz)<\fz$, we then have
$$\begin{array}{cl}
\wz(\{x\in\rz: |[b,T]f(x)|>\lz\})&= \wz(\{x\in\rz:
|\dsum_{j\in\nn_0}\lz_j[b,T]a_j(x)|>\lz\})\\
&\le \wz(\{x\in\rz:
|\dsum_{j\in E_1}\lz_j[b,T]a_j(x)|>\lz/3\})\\
&\quad+\wz(\{x\in\rz:
|\dsum_{j\in E_2}\lz_j[b,T]a_j(x)|>\lz/3\})\\
&\quad+\wz(\{x\in\rz:
|\lz_0[b,T]a_0(x)|>\lz/3\})\\
&:=F_1+F_2+F_3,
\end{array}$$
where $E_1=\{j\in\nn:\ r_j<\dz_0\}$ and $E_2=\{j\in\nn:\ 2\ge r_j\ge
\dz_0\}$ and $\dz_0$ be a number satisfying
$4\dz_0=\dz_0^{1/(1+\tz)}$. Obviously, $\dz_0<1$.

For $F_1$, let $b_j=\frac1{|Q_j|}\int_{Q_j}b(y)dy$. Note that
$$\begin{array}{cl}
\dsum_{j\in E_1}\lz_j[b,T]a_j(x)&=\dsum_{j\in
E_1}\lz_j[b-b_j,T]a_j(x)\\
&=\dsum_{j\in
E_1}\lz_j(b(x)-b_j)Ta_j(x)\chi_{4nQ_j}(x)\\
&\quad+\dsum_{j\in
E_1}\lz_j(b(x)-b_j)Ta_j(x)\chi_{(4nQ_j)^c}(x)\\
&\quad-T(\dsum_{j\in E_1}\lz_j(b(x)-b_j)a_j)(x)\\
&:=F_{11}(x)+F_{12}(x)+F_{13}(x).
\end{array}$$
Thus, by (i) of Theorem 7.2 and Theorem 7.1, we obtain
$$\begin{array}{cl}
\wz(\{x\in\rz:\ |F_{11}(x)|>\lz/9\})&\le \dfrac C\lz\dsum_{j\in
E_1}|\lz_j|\|(b-b_j)(Ta_j)\chi_{4nQ_j}\|_{L^1_\wz(\rz)}\\
&\le \dfrac C\lz\dsum_{j\in
E_1}|\lz_j|\|(b-b_j)\chi_{4nQ_j}\|_{L^2_\wz(\rz)}\|a_j\|_{L^2_\wz(\rz)}\\
&\le \dfrac C\lz\dsum_{j\in\nn}|\lz_j|\|b\|_{BMO^{loc}}\\
&\le \dfrac C\lz\|b\|_{BMO^{loc}}\|f\|_{h^1_\wz(\rz)}.
\end{array}$$
By the weighted weak type (1,1) of $T$ (see Theorem 7.1 (ii)) and Proposition 7.1, we
get
$$\begin{array}{cl}
\wz(\{x\in\rz:\ |F_{13}(x)|>\lz/9\})&\le \dfrac C\lz\dsum_{j\in
E_1}|\lz_j|\|(b-b_j)a_j\|_{L^1_\wz(\rz)}\\
&\le \dfrac C\lz\dsum_{j\in \nn}|\lz_j|\|b\|_{BMO^{loc}}\\
&\le \dfrac C\lz\|b\|_{BMO^{loc}}\|f\|_{h^1_\wz(\rz)}.
\end{array}$$
Now we consider the term $F_{12}(x)$. Obviously,
$$
\wz(\{x\in\rz:\ |F_{12}(x)|>\lz/6\}|\le \lz^{-1}\dsum_{j\in
E_1}|\lz_j|\dint_\rz T((b(x)-b_j)a_j)(x)\wz(x)dx.\eqno(7.4)$$ We
claim that
$$\dint_\rz T((b(x)-b_j)a_j)(x)\wz(x)dx\le C$$
holds for all atoms $a_j$ for $j\in E_1$. For convenience, we
denote  $a_j$ by $a$, $Q_j(x_j,r_j)$ by $Q(x_0,\dz)$ and $b_j$ by
$b_Q$ for $j\in E_1$ .

 We let $Q^*=4Q$
and $\bar Q=Q(x_0,\dz^{1/(1+\tz)})$. Then
$$\begin{array}{cl}
\dint_\rz|T(b-b) a|\wz(x)dx&\le \dint_{Q^*}|T(b-b_Q)
a|\wz(x)dx\\
&\qquad+ \dint_{\bar Q\setminus Q^*}|T(b-b_Q)
a|\wz(x)dx+\dint_{\rz\setminus \bar Q}|T
a|\wz(x)dx\\
&:=I+II+III.\end{array}$$ For $I$, similar to (7.2), we have
$$I\le
C\l(\dint_{\rz}|T(b-b_Q)a|^p\wz(x)dx\r)^{1/p}\l(\dint_{Q^*}\wz(x)dx\r)^{1/p'}\le
C.$$ We now estimate the term $III$. Clearly, by the mean value
theorem,
$$\begin{array}{cl}
|T(b-b_Q)a(x)|&\le \dfrac {C\dz|b(x)-b_Q|}{|x-x_0|^{\tz+n+1}}\chi_{\{|x-x_0|<4n\}}(x)\dint_Q|a(y)|dy\\
&\le \dfrac {C\dz|b(x)-b_Q|}{|x-x_0|^{\tz+n+1}}\chi_{\{|x-x_0|<4n\}}(x)\\
&\qquad\times\l(\dint_Q|a(x)|^p\wz(x)dy\r)^{1/p}\l(\dint_{Q}[\wz(x)]^{-p'/p}dx\r)^{1/p'}\\
&\le \dfrac
{C\dz|b(x)-b_Q|}{|x-x_0|^{\tz+n+1}}\chi_{\{|x-x_0|<4n\}}(x)\dfrac
{|Q|}{\wz(Q)}.
\end{array}$$
Hence, by the properties of $A_1^{loc}$ (see Lemma 2.1), we have
$$\begin{array}{cl}
III&\le \dfrac
{C\dz|Q|}{\wz(Q)}\dint_{\dz^{1/(1+\tz)}\le|x-x_0|\le
4n}\dfrac {|b(x)-b_Q|\wz(x)}{|x-x_0|^{\tz+n+1}}dx\\
&\le \dfrac {C\dz|Q|}{\wz(Q)}\dsum_{k=k_0}^{k_1}\dfrac
1{(2^k\dz)^{1+\tz}}\l(\dfrac 1{(2^k\dz)^n}\dint_{|x-x_0|\le
2^k\dz}\dfrac {|b(x)-b_Q|\wz(x)}{|x-x_0|^{\tz+n+1}}dx\r)\\
&\le C,
\end{array}$$
where $k_0$  and $k_1$ are positive integers such that
$2^{k_0}\dz\le \dz^{1/(1+\tz)}\le 2^{k_0+1}\dz$ and $2^{k_1-1}\le
4n\le 2^{k_1}$. We now estimate the term $II$. For $x\in \bar
Q\setminus Q^*$
$$\begin{array}{cl}
Ta(x) &=\dint_\rz\dfrac
{e^{i|x-y|^{-\tz}}v(x-y)}{|x-y|^{n(2+\tz)/r'}}\\
&\qquad\times\l(\dfrac1{|x-y|^{n(1-(2+\tz)/r')}}
-\dfrac1{|x_0-x|^{n(1-(2+\tz)/r')}}\r)(b(x)-b_Q)a(y)dy\\
&\qquad+\dint_\rz\dfrac
{e^{i|x-y|^{-\tz}}v(x-y)}{|x-y|^{n(2+\tz)/r'}}\dfrac{(b(x)-b_Q)a(y)}{|x_0-x|^{n(1-(2+\tz)/r')}}dy\\
&=A(x)+B(x),
\end{array}$$
where $r'$ is taken so close to 1 to guarantee that $2+\tz<r$.
Applying the mean value theorem to the term in brackets in the
integrand of $A$, and noting that for $y\in Q$, and $x\in \bar
Q\setminus Q^*$, $|x-y|\ge c|x-x_0|$, we have
$$\begin{array}{cl}
|A(x)|&\le \dfrac
{C\dz|b(x)-b_Q|}{|x-x_0|^{n+1}}\chi_{\{|x-x_0|<4n\}}(x)\dint_Q|a(y)|dy\\
&\le \dfrac
{C\dz|b(x)-b_Q|}{|x-x_0|^{n+1}}\chi_{\{|x-x_0|<4n\}}(x)\dfrac
{|Q|}{\wz(Q)}.\end{array}$$ Therefore,
$$\begin{array}{cl}
II&\le \dfrac {C\dz|Q|}{\wz(Q)} \dint_{\dz\le|x-x_0|\le
4n}\dfrac {|b(x)-b_Q|\wz(x)}{|x-x_0|^{n+1}}dx\\
&\qquad +C\dint_{\dz\le|x-x_0|<\dz^{1/(1+b)}}|K_{\tz,r}*a|\dfrac{|b(x)-b_Q|\wz(x)}{|x_0-x|^{n(1-(2+\tz)/r')}}dx\\
&\le C+C\l(\dint_\rz
|K_{\tz,r}*a|^rdx\r)^{1/r}\l(\dint_{\dz<|x-x_0|<\dz^{1/(1+\tz)}}
\dfrac{|b(x)-b_Q|\wz(x)^{r'}}{|x_0-x|^{n(1-(2+\tz)/r')}}dx\r)^{1/r'}\\
&\le C.
\end{array}$$
Thus, the claim is proved.

From (7.4), we have
$$\wz(\{x\in\rz: |F_{12}(x)|>\lz/9\})\le C\lz^{-1}\dsum_{j\in E_1}|\lz_j|
\le C\lz^{-1}\|f\|_{h^1_\wz(\rz)}.$$ It remains to consider the
term $F_2$. In fact, it is very simple. Let $Q^*_j=\frac
{10n}{\dz_0} Q_j$. Note that for any atom $a_j$
$$\dint_\rz |[b,T]a_j|\wz(x)dx=\dint_{Q^*_j} |[b,T]a_j|\wz(x)dx.$$ Note that $\wz(Q_j^*)\le C\wz(Q_j)$, we
then have
$$\begin{array}{cl}
F_2&=\wz(\{x\in\rz:
|\dsum_{j\in E_2}\lz_j[b,T]a_j(x)|>\lz/3\})\\
&\le {\lz}^{-1}\dsum_{j\in E_2}|\lz_j|\dint_\rz |[b,T]a_j|\wz(x)dx\\
&\le {\lz}^{-1}\dsum_{j\in E_2}|\lz_j|\dint_{Q^*_j} |[b,T]a_j|\wz(x)dx\\
&\le {\lz}^{-1}\dsum_{j\in E_2}|\lz_j|\|[b,T]a_j\|_{L^2_\wz(\rz)}[\wz(Q^*_j)]^{1/2}\\
&\le C\|b\|_{BMO^{loc}}{\lz}^{-1}\dsum_{j\in E_2}|\lz_j|\|a_j\|_{L^2_\wz(\rz)}[\wz(Q^*_j)]^{1/2}\\
&\le C\|b\|_{BMO^{loc}}{\lz}^{-1}\dsum_{j\in E_2}|\lz_j|.
\end{array}$$
It remains to estimate the term $F_3$.
$$\begin{array}{cl}F_3&\le C\dfrac {|\lz_0|}\lz\dint_\rz|[b,T]a(x)|\wz(x)dz\\
&\le C\dfrac {|\lz_0|}\lz\|[b,T]a\|_{L^2_\wz(\rz)}[\wz(\rz)]^{1/2}\\
&\le C\dfrac {|\lz_0|}\lz\|a\|_{L^2_\wz(\rz)}[\wz(\rz)]^{1/2}\le
C\dfrac {|\lz_0|}\lz.
\end{array}$$
 From these, we have
$$\begin{array}{cl}
\wz(\{x\in\rz: |[b,T]f(x)|>\lz\})&\le  \dsum_{i=1}^3|\{x\in\rz:\
|F_{1i}(x)|>\lz/9\}|\\
&\qquad+ \wz(\{x\in\rz:\
|F_{2}(x)|>\lz/3\})\\
&\qquad+ \wz(\{x\in\rz:\
|F_{3}(x)|>\lz/3\})\\
&\le \dfrac C\lz\|b\|_{BMO^{loc}}\|f\|_{h^1_\wz(\rz)}.
\end{array}$$
Thus, the proof of Theorem 7.2 is complete.

Next we  show that the pseudodifferential operators are bounded on
$h^p_\wz(\rz)$, where the weight $\wz$ is in the weight class
$A_p(\vz)$  which  is  contained in $A_p^{\loc}$ for $1\le p<\fz$.
Let us first introduce some definitions.

Let $m$ be real number. Following \cite{t}, a symbol in
$S_{1,\dz}^m$ is a smooth function $\sz(x,\xi)$ defined on
$\rz\times\rz$ such that for all multi-indices $\az$ and $\bz$ the
following estimate holds:
$$|D_x^\az D_\xi^\bz\sz(x,\xi)|\le
C_{\az,\bz}(1+|\xi|)^{m-|\bz|+\dz|\az|},$$ where $ C_{\az,\bz}>0$
is independent of $x$ and $\xi$.

The operator $T$ given by
$$Tf(x)=\dint_\rz\sz(x,\xi)e^{2\pi ix\cdot \xi}\hat f(\xi)\,d\xi$$
is called a pseudo-differential operator with symbol
$\sz(x,\xi)\in S_{1,\dz}^m$, where $f$ is a Schwartz function and
$\hat f$ denotes the Fourier transform of $f$.

 In the rest of this section, we let
$\vz(t)=(1+t)^\az$ with $\az>0$.

 A weight will always mean a positive function which is
locally integrable. We say that a weight $\wz$ belongs to the class
$A_p(\vz)$ for $1<p<\fz$, if there is a constant $C$ such that for
all cubes  $Q=Q(x,r)$ with center $x$ and sidelength $r$
$$\l(\dfrac 1{\vz(|Q|)|Q|}\dint_Q\wz(y)\,dy\r)
\l(\dfrac 1{\vz(|Q|)|Q|}\dint_Q\wz^{-\frac 1{p
-1}}(y)\,dy\r)^{p-1}\le C.$$
 We also
say that a  nonnegative function $\wz$ satisfies the $A_1(\vz)$
condition if there exists a constant $C$ for all cubes $Q$
$$M_\vz(\wz)(x)\le C \wz(x), \ a.e.\ x\in\rz.$$
 where
$$M_\vz f(x)=\dsup_{x\in Q}\dfrac 1{\vz(|Q|)|Q|}\dint_Q|f(y)|\,dy.$$
   Since $\vz(|Q|)\ge 1$, so
$A_p(\rz)\subset A_p(\vz)$ for $1\le p<\fz$, where $A_p(\rz)$
denote the classical Muckenhoupt weights; see \cite{gr}.

{\bf Remark :}\quad It is easy to see that if $\wz\in A_p(\vz)$,
then $\wz(x)dx$ may be not a doubling measure. In fact, let
$\az>0$ and $0\le\gz<\az$, it is easy to check that
$\wz(x)=(1+|x|\log(1+|x|))^{-(n+\gz)}\not\in A_\fz(\rz)$ and
$\wz(x)dx$ is not a doubling measure, but
$\wz(x)=(1+|x|\log(1+|x|))^{-(n+\gz)}\in A_1(\vz)$ provided that
$\vz(r)=(1+r^{1/n})^\az$.

Similar to the classical Muckenhoupt weights, we give some
properties for weights $\wz\in A_\fz(\vz)=\bigcup_{p\ge
1}A_p(\vz)$.
\begin{lem}\label{l7.3.}\hspace{-0.1cm}{\rm\bf 7.3.}\quad For any cube
$Q\subset \rz$, then
\begin{enumerate}
\item[(i)]If $ 1\le p_1<p_2<\fz$, then $A_{p_1}(\vz)\subset
A_{p_2}(\vz)$. \item[(ii)] $\wz\in A_p(\vz)$ if and only if
$\wz^{-\frac 1{p-1}}\in A_{p'}(\vz)$, where $1/p+1/p'=1.$
\item[(iii)] If $\wz\in A_p$ for $1\le p<\fz$, then for any
measurable set $E\subset Q$,
$$\dfrac {|E|}{\vz(|Q|)|Q|}\le C\l(\dfrac
{\wz(E)}{\wz(Q)}\r)^{1/p}.$$
\end{enumerate}
\end{lem}

\begin{lem}\label{l7.4.}\hspace{-0.1cm}{\rm\bf 7.4.}\quad
Let $T$ be the $S_{1,0}^0$ pseudodifferential operators, then
$$\|Tf\|_{L^p_\wz(\rz)}\le C_{p,\wz}\|f\|_{L^p_\wz(\rz)}$$ for
$1<p<\fz$ and $\wz\in A_p(\vz)$.
\end{lem}
 Lemmas 7.3 and   7.4 can be founded in \cite{t}. The following lemma was
proved in \cite{g}.
\begin{lem}\label{l7.5.}\hspace{-0.1cm}{\rm\bf 7.5.}\quad
Let $T$ be the $S_{1,0}^0$ pseudodifferential operators, if
$\vz\in {\cal D}$ then $T_tf=\vz_t*Tf$ has a symbol $\sz_t$ which
satisfies $D^\bz_xD_\xi^\az\sz_t(x,\xi)\le
C_{\az,\bz}(1+|\xi|)^{-|\az|}$ and a kernel
$K_t(x,z)=FT_\xi\sz_t(x,\xi)$ which satisfies $|D^\bz_x
D^\az_z K_t(x,z)|\le C_{\az,\bz}|z|^{-n-|\az|}$, where $C_{\az,\bz}$
is independent of $t$ if $0<t<1$.
\end{lem}

\begin{thm}\label{t7.3.}\hspace{-0.1cm}{\rm\bf 7.3.}\quad
Let   $T$ be the $S_{1,0}^0$ pseudodifferential
operators, then
$$\|Tf\|_{h^p_\wz(\rz)}\le C_{p,\wz}\|f\|_{h^p_\wz(\rz)}$$ for $\wz\in
A_\fz(\vz)$ and $0<p\le 1$.
\end{thm}
Proof. Since  $\wz\in A_\fz(\vz)$, so  $\wz\in A_q(\vz)$ for some
$q>1$. By Theorem  6.2, it suffices to show that for any atom
$(p,q,s)_\wz$ $a$ supported $Q=Q(x_0,r)$ with $r\le 2$ and
$\|a\|_{L^q_\wz(\rz)}\le [\wz( Q)]^{1/p-1/q}$, such that
$$\|Ta\|_{h^p_\wz(\rz)}\le C_{\wz,p},\eqno(7.6)$$  and if $a$ is a
single atom, then
$$\|Ta\|_{h^p_\wz(\rz)}\le C_{\wz,p}.\eqno(7.7)$$
Obviously, (7.7) holds. Now we prove (7.6).

 If $Q^*=2Q$, we then have
$$\begin{array}{cl}
\dint_{Q^*}\dsup_{0<t<1}|\vz_t*Ta(x)|^p\wz(x)dx&\le
\wz(Q^*)^{(q-p)/q}\l(\dint_{Q^*}\dsup_{0<t<1}|\vz_t*Ta(x)|^q\wz(x)dx\r)^{p/q}\\
&\le C\wz(Q^*)^{(q-p)/q}\l(\dint_{\rz}|Ta|^q\wz(x)dx\r)^{p/q}\\
&\le C\wz(Q^*)^{(q-p)/q}\l(\dint_{\rz}|a|^q\wz(x)dx\r)^{p/q}\\
&\le C.
\end{array}$$
 To estimate $\int_{\rz\setminus
Q^*}\dsup_{t<1}|\vz*Ta|^p$, we consider two cases.

The first case is when $r<1$. We expand $K_t(x,x-z)$ in a Taylor
series about $z=x_0$ so that
$$\vz_t*Ta(x)=\dint_\rz
K_t(x,x-z)a(z)dx=\dint_\rz\dsum_{\|\az|=N+1}D_z^\az
K_t(x,x-\xi)z^\az a(z)dz,$$ where $\xi$ is in $Q$,  and hence by Lemma 7.5,
$$|\vz_t*Ta(x)|\le C|x-x_0|^{-(n+N+1)}\|a\|_1 |Q|^{(N+1)/n}.$$ Taking
$N$ is large enough and $r<1$, by Lemma 7.3 (iii), we  then have
$$\begin{array}{cl}
\dint_{\rz\setminus
Q^*}&\dsup_{0<t<1}|\vz_t*Ta(x)|^p\wz(x)dx\\
&\le
C|Q|^{p(N+1)/n}\dfrac{|Q|^p}{\wz(Q)}\dint_{\rz\setminus
Q^*}|x-x_0|^{-p(n+N+1)}\wz(x)dx\\
&\le C|Q|^{p(N+1)/n}\dfrac{|Q|^p}{\wz(Q)} \dsum_{k=1}^\fz
(2^kr)^{-p(n+N+1)}\dint_{|x-x_0|<2^kr}\wz(x)dx\\
\end{array}$$
$$\begin{array}{cl}
&\le C\dfrac 1{\wz(Q)} \dsum_{k=1}^{k_0}(2^k)^{-p(n+N+1)}\dint_{|x-x_0|<2^kr}\wz(x)dx\\
&\quad+C\dfrac1{\wz(Q)} \dsum_{k=k_0}^\fz (2^k)^{-(n+N+1)}\dint_{|x-x_0|<2^kr}\wz(x)dx\\
&\le C\dfrac1{\wz(Q)} \dsum_{k=1}^{k_0}2^{knq}2^{-kp(n+N+1)}\wz(Q)\\
&\qquad
+C\dfrac1{\wz(Q)} \dsum_{k=k_0}^\fz (2^kr)^{-p(n+N+1)+\az n}2^{knq}\wz(Q)\\
 &\le C,\end{array}$$ where the integer $k_0$ satisfies  $2^{k_0-1}\le 1/r<2^{k_0}$. To estimate
with the case when $1<r\le2$, by Lemma 7.5, for all $M>0$, we have
$$|K_t(x,x-z)|\le C_M|x-z|^{-M}.$$ So
$$|\vz_t*Ta(x)|\le\dint_Q
|K_t(x,x-z)a(z)|dz\le C_M|x-x_0|^{-M}\|a\|_{L^1(\rz)}.$$  Note
that $1<r\le 2$, we then have
$$\begin{array}{cl}
\dint_{\rz\setminus Q^*}\dsup_{0<t<1}|\vz_t*Ta(x)|^p\wz(x)dx &\le
C_M\|a\|^p_{L^1(\rz)}\dint_{\rz\setminus Q^*}|x-x_0|^{-Mp}\wz(x)dx\\
&\le C_M\dfrac{|Q|^p}{\wz(Q)}\dint_{\rz\setminus Q^*}|x-x_0|^{-M}\wz(x)dx \\
&\le C_M\dfrac1{\wz(Q)} \dsum_{k=0}^\fz (2^kr)^{-Mp+\az
n}2^{knp}\wz(Q)\\
&\le C,\end{array}$$ if $M$ is large enough. The proof is
complete.

 LMAM, School of Mathematical  Science

 Peking University

 Beijing, 100871

 P. R. China

\bigskip

 E-mail address:  tanglin@math.pku.edu.cn

\end{document}